\documentclass[12pt]{article}
\usepackage{epsfig}
\usepackage{amsmath}
\usepackage{amssymb}
\usepackage{eucal}

\newtheorem{thm}{Theorem}[section]

\newtheorem{hyp}[thm]{Hypothesis}

\newcommand{\Section}[1]{\section{#1} \setcounter{equation}{0}}

\title{A degenerate Hopf bifurcation in retarded functional differential equations, and applications to endemic bubbles.}

 \author{
Victor G. LeBlanc\\Department of
  Mathematics and Statistics\\University of Ottawa\\Ottawa, 
ON K1N 6N5\\CANADA}
\date{\today}

\begin{document}

\maketitle

\begin{abstract}
In this paper, we study degenerate Hopf bifurcations in a class of parametrized retarded functional differential equations.  Specifically, we are interested in the case where the eigenvalue crossing condition of the classical Hopf bifurcation theorem is violated.  Our approach is based on center manifold reduction and Poincar\'e-Birkhoff normal forms, and a singularity theoretical classification of this degenerate Hopf bifurcation.  Our results are applied to a recently developed SIS model incorporating a delayed behavioral response.  We show that the phenomenon of {\it endemic bubbles}, which is characterized by a branch of periodic solutions which bifurcates from the endemic equilibrium at some value of the basic reproduction number $R_0$, and then reconnects to the endemic equilibrium at a larger value of $R_0$, originates in a codimension-two organizing center where the eigenvalue crossing condition for the Hopf bifurcation theorem is violated.
\end{abstract}

\pagebreak
\Section{Introduction}

Retarded functional differential equations (RFDEs), of which delay differential equations (DDEs) are a special case, are used to model a large variety of phenomena in sciences, engineering, economics and many other areas  
\cite{BBL,BR,BPS,HFEKGG,Kuang,LM,MT,NOC,PPK,SieberKrauskopf,SC,SS,VTK}.
One of the main technical differences between RFDEs and ordinary differential equations (ODEs) is that while ODEs require initial data only at one point (typically at time $t=0$) to generate a solution, RFDEs require initial data in a range of past values, typically an interval $t\in [-\tau,0]$.  One of the consequences of this fact is that the phase space for RFDEs is infinite-dimensional.  Despite this fact, most of the usual tools and techniques of the theory of dynamical systems can be suitably adapted to the study of RFDEs \cite{FM1,FM2,HVL}.  In particular, for parametrized RFDEs, one can analyze bifurcations using center manifold theory and normal forms \cite{FM1,FM2}, and in many studies (see for example \cite{BR,BPS,MT,NOC,SieberKrauskopf,SC}), this has led to valuable insight into many phenomena which are modeld using RFDEs.

The focus of this paper is an analysis of a certain codimension 2 degenerate Hopf bifurcation in RFDEs.  In particular, we are interested in studying the case where a parametrized RFDE admits an equilibrium solution which satisfies the simple purely imaginary eigenvalue condition of the Hopf bifurcation theorem \cite{MarsMC}, but for which the crossing condition of this theorem is violated.  This is motivated from recent results \cite{bubbles} where a SIS model incorporating a delayed behavioral response was analyzed and shown to exhibit a phenonmenon which the authors called {\it endemic bubbles}.  It became apparent to us that the bifurcation diagrams reported in \cite{bubbles} resembled bifurcation diagrams which had been found to be in the versal unfolding of the degenerate Hopf bifurcation (with crossing condition violated) in \cite{GoLang}.  Although this SIS model is the application which motivated our analysis, the theoretical results we present here have a much larger scope of application, since many phenomena in nature are modelled using the class of delay differential equations we study.

\subsection{The crossing condition for the Hopf bifurcation}

Consider as a prototype the delay differential equation
\begin{equation}
\dot{x}(t)=\alpha\,x(t)+\beta\,x(t-\tau)+F(x(t),x(t-\tau))
\label{lindde1}
\end{equation}
where $\alpha$, $\beta\in\mathbb{R}$ are parameters, $\tau>0$ is a fixed delay time, and $F$ represents higher-order nonlinear terms.  

Equation (\ref{lindde1}) has $x=0$ as an equilibrium solution for all values of $\alpha$ and $\beta$.  We are interested in bifurcations from this trivial equilibrium point.  To this end, we consider the characteristic equation
\[
\xi=\alpha+\beta\,e^{-\xi\tau}.
\]
In particular, we are interested in purely imaginary solutions $\xi=i\,\omega$ to this characteristic equation.  It is easy to see that such solutions occur when the parameters $\alpha$ and $\beta$ satisfy 
\begin{equation}
\alpha+\beta\cos\,\tau\omega=0,\,\,\,\,\,-\omega=\beta\sin\,\tau\omega,\,\,\,\,\,\omega=\sqrt{\beta^2-\alpha^2},\,\,\,\,\,\beta^2>\alpha^2.
\label{bifcurve}
\end{equation}
\begin{figure}[htpb]
\begin{center}
\includegraphics[width=2.5in]{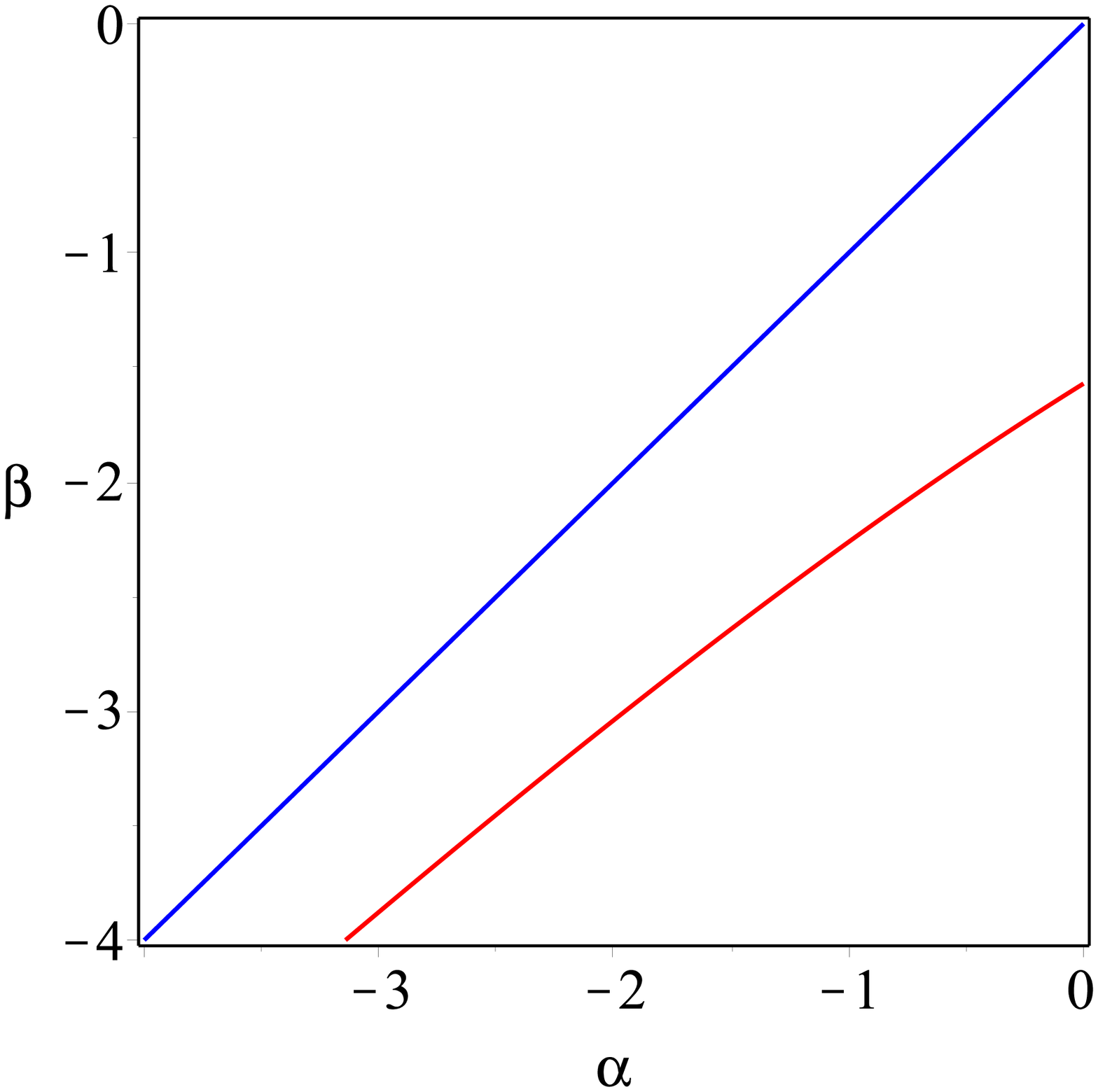}
\includegraphics[width=2.5in]{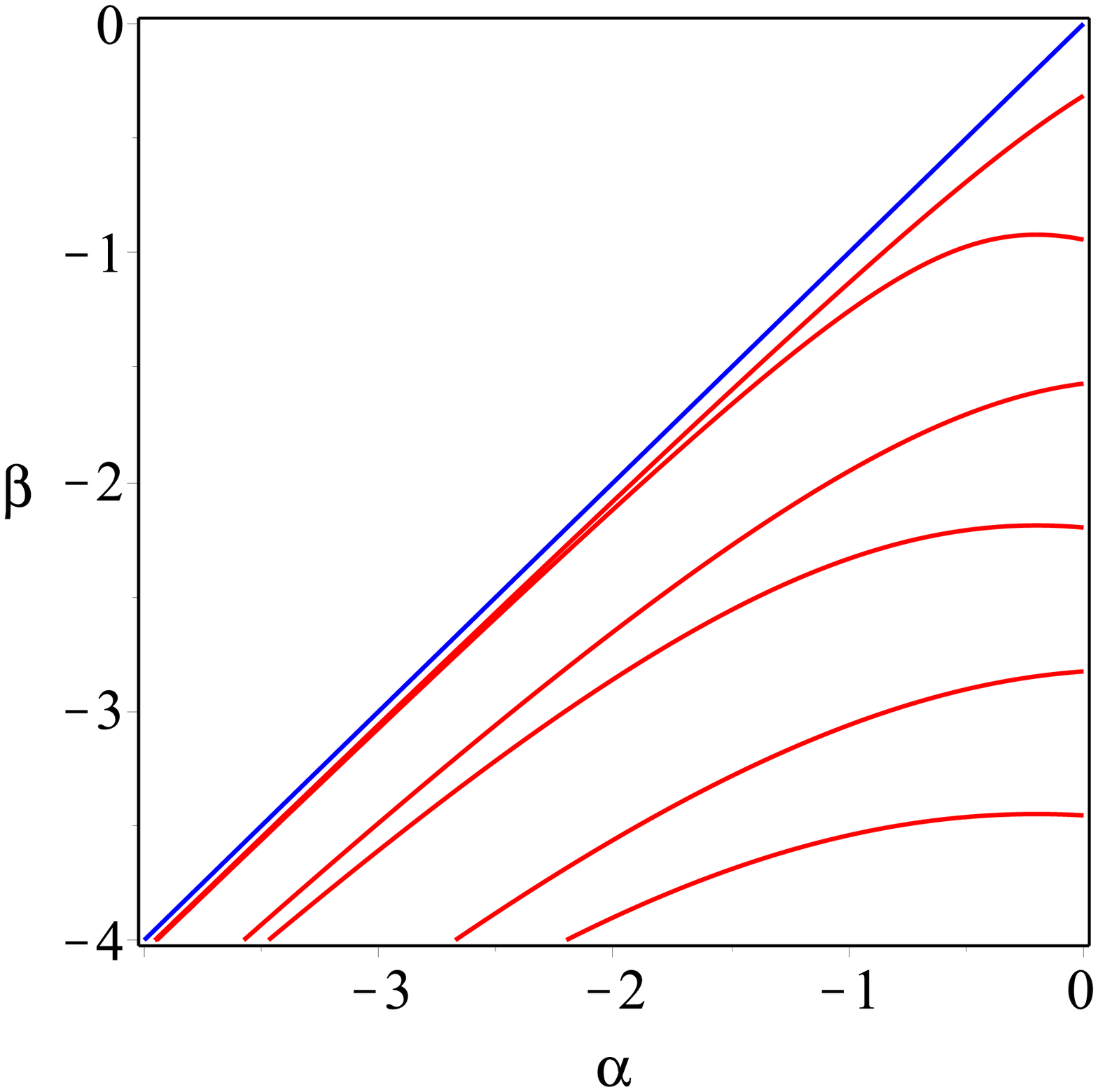}
\caption{Hopf bifurcation curves (\ref{bifcurve}) in red for $\tau=1$ on the left and $\tau=5$ on the right.  The blue line represents $\beta=\alpha$.  In the region in parameter space below the blue line and above the uppermost red curve, the equilibrium point $x=0$ of (\ref{lindde1}) is locally asymptotically stable, and loses this stability below the uppermost red curve.}
\label{bifcurve1}
\end{center}
\end{figure}
Equations (\ref{bifcurve}) define curves in the $\alpha$-$\beta$ parameter space, as illustrated in Figure \ref{bifcurve1}.  

\begin{figure}[htpb]
\begin{center}
\includegraphics[width=2.5in]{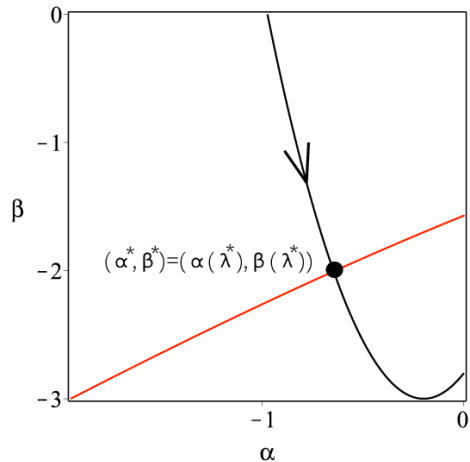}
\caption{As the control parameter $\lambda$ varies through $\lambda^*$, the (black) curve $(\alpha(\lambda),\beta(\lambda))$ crosses the (red) Hopf bifurcation curve (\ref{bifcurve}) at the point $(\alpha^*,\beta^*)$.  Generically, the crossing is transversal (as illustrated here), and we get a Hopf bifurcation from the trivial equilibrium point of (\ref{lindde1}).}
\label{fig2}
\end{center}
\end{figure}

Now let us suppose that $\alpha$ and $\beta$ in (\ref{lindde1}) depend on a distinguished external control parameter $\lambda$.   If the path $(\alpha(\lambda),\beta(\lambda))$ in parameter space crosses the Hopf bifurcation curve (\ref{bifcurve}) at the point $(\alpha^*,\beta^*)=(\alpha(\lambda^*),\beta(\lambda^*))$ from the region where the equilibrium $x=0$ is stable into the region where it is unstable as illustrated in Figure \ref{fig2}, then typically it crosses this Hopf curve transversally, and assuming non-degeneracy conditions in the nonlinear coefficients (specifically, the first Lyapunov coefficient is non-zero), a super- or sub-critical Hopf bifurcation occurs from the trivial equilibrium.

\begin{figure}[htpb]
\begin{center}
\includegraphics[width=2.5in]{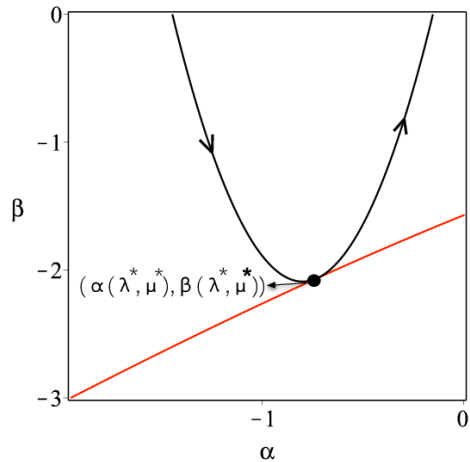}
\caption{Tangential intersection of the parameter curve $(\alpha(\lambda,\mu^*),\beta(\lambda,\mu^*))$ with the Hopf bifurcation curve (\ref{bifcurve}) at $\lambda=\lambda^*$.  The crossing condition of the Hopf bifurcation theorem is violated for this parameter path.}
\label{fig3}
\end{center}
\end{figure}

\begin{figure}[htpb]
\begin{center}
\includegraphics[width=2in]{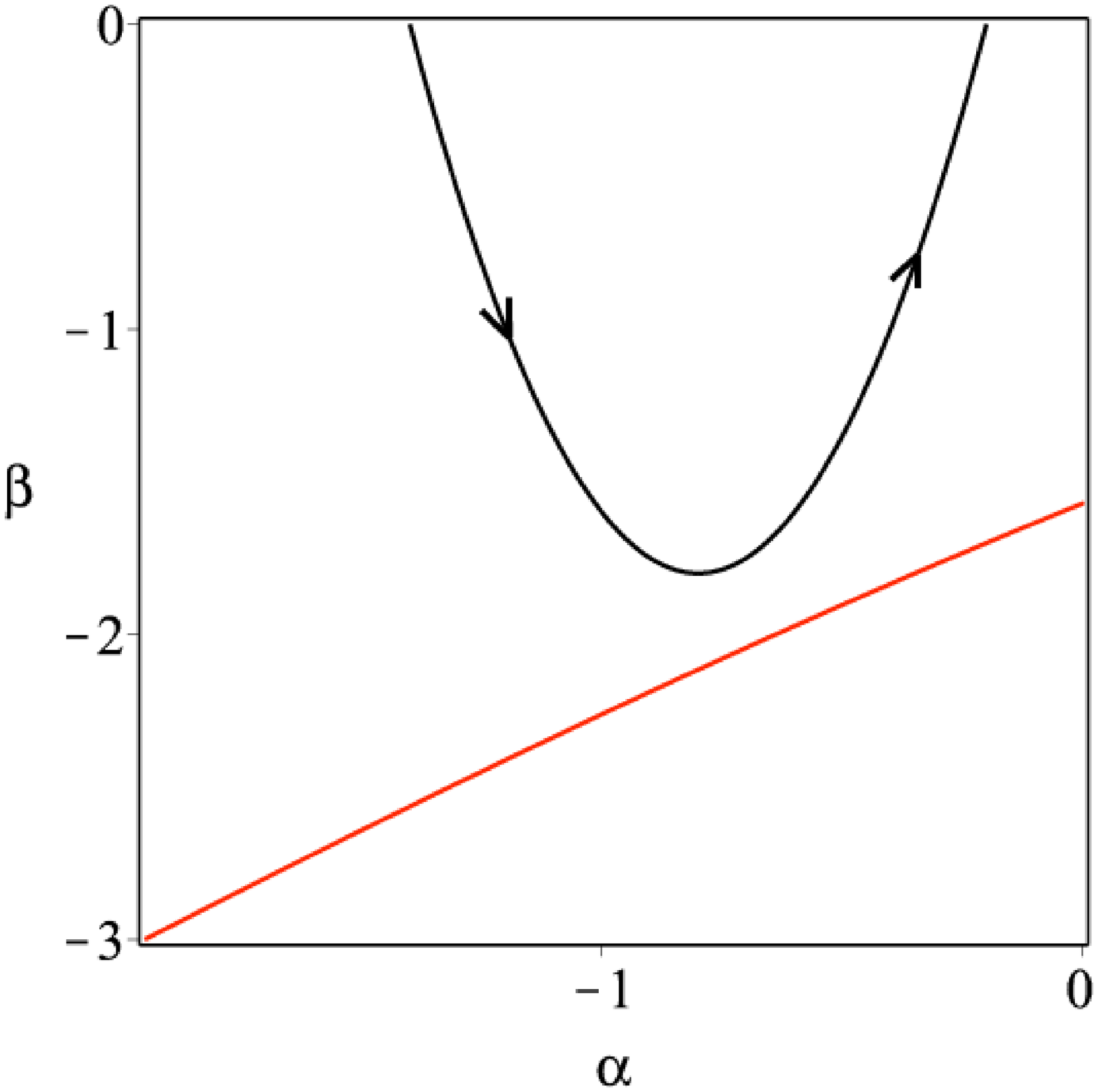}
\includegraphics[width=2in]{DEGEN2.eps}
\includegraphics[width=2in]{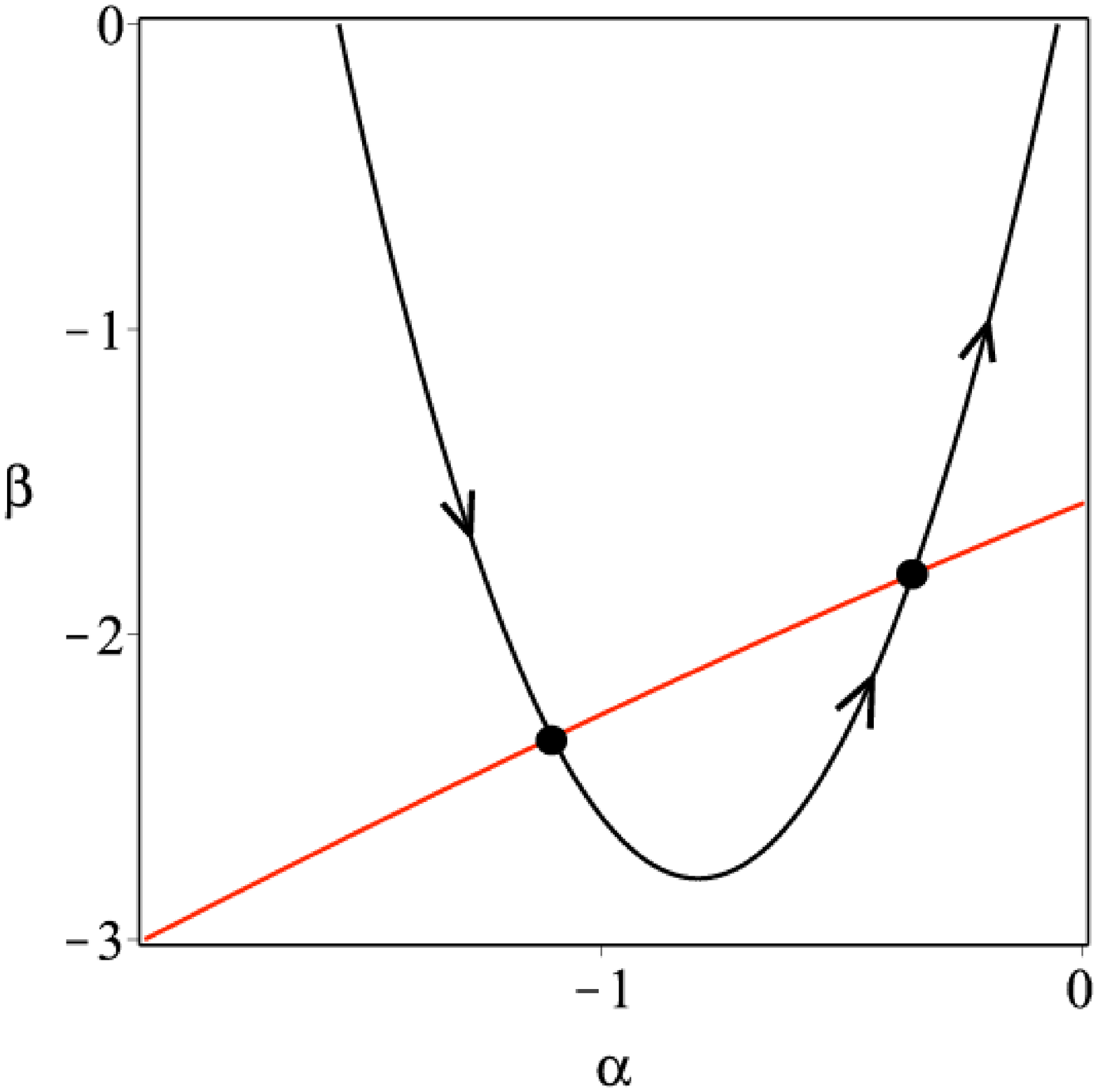}
\caption{Unfolding the tangential intersection (middle diagram), for fixed values of $\mu$ in a neighborhood of $\mu^*$.  Generically, we either get no intersections (as illustrated on the left) or two intersections (as illustrated on the right) with the Hopf bifurcation curve.}
\label{fig4}
\end{center}
\end{figure}

\begin{figure}[htpb]
\begin{center}
\includegraphics[width=5.5in]{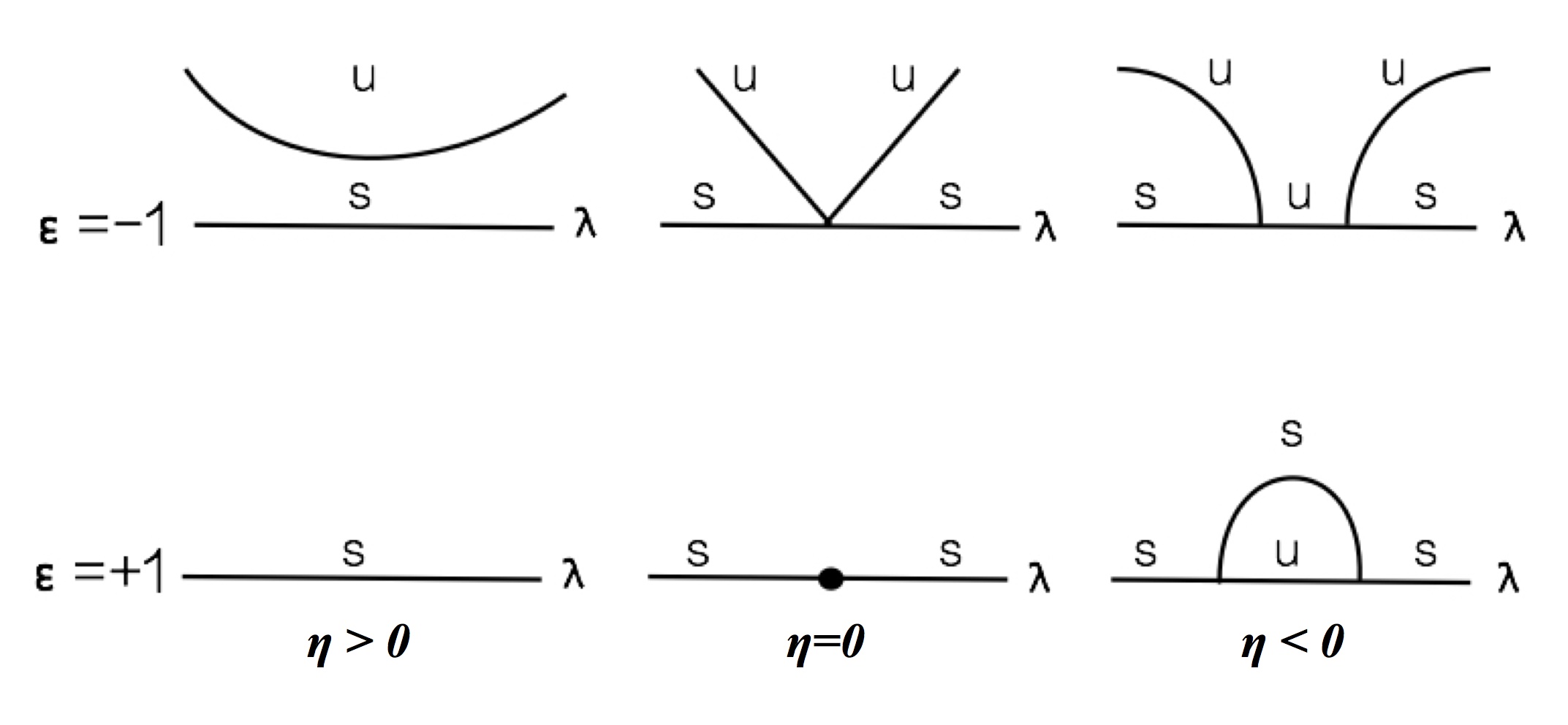}
\caption{Universal unfolding of the degenerate crossing condition for the Hopf bifurcation.  These are the bifurcation diagrams for the normal form (\ref{normform}).  An ``s'' designates a stable branch of periodic solutions, whereas a ``u'' designates an unstable branch of periodic solutions.}
\label{fig5}
\end{center}
\end{figure}

In applications, it often occurs that in addition to depending on the distinguished parameter $\lambda$, the coefficients $\alpha$ and $\beta$ in (\ref{lindde1}) also depend on auxiliary parameters $\mu\in\mathbb{R}^p$.  For the purposes of this discussion, let us suppose that $p=1$.
In this case, it is possible that for  a certain value of $\mu=\mu^*$, the curve $\lambda\longmapsto (\alpha(\lambda,\mu^*),\beta(\lambda,\mu^*))$ intersects the Hopf curve (\ref{bifcurve}) tangentially at the point $(\alpha^*,\beta^*)=(\alpha(\lambda^*,\mu^*),\beta(\lambda^*,\mu^*))$ as is illustrated in Figure \ref{fig3}.  In this case, the transversality condition of the Hopf Bifurcation Theorem \cite{MarsMC} is violated.  However, as $\mu$ varies in a neighborhood of $\mu^*$, this degeneracy is ``unfolded'' as illustrated in Figure \ref{fig4}.  Assuming certain non-degeneracy conditions on nonlinear terms for ODEs undergoing a violation of the crossing condition, it is shown in \cite{GoLang} that the Hopf bifurcation diagrams in a neighborhood of $\mu^*$ are equivalent to one of the bifurcation diagrams of the normal form
\begin{equation}
x(\varepsilon(\lambda^2+\eta)+x^2)=0,\,\,\,\,\,\,\varepsilon=\pm 1
\label{normform}
\end{equation}
as shown in Figure \ref{fig5}.


Consider the twice differentiable path $\lambda\longmapsto (\alpha(\lambda,\mu),\beta(\lambda,\mu))$, and let $\lambda^*$ and $\mu^*$ be such that the point
\[
(\alpha^*,\beta^*)\equiv (\alpha(\lambda^*,\mu^*),\beta(\lambda^*,\mu^*))
\]
satisfies both the following equations
\begin{eqnarray}
{\displaystyle\alpha^*+\beta^*\cos\,\tau\sqrt{\beta^{*2}-\alpha^{*2}}}&=&0\label{cond1}\\[0.1in]
{\displaystyle\beta^*\alpha_{\lambda}(\lambda^*,\mu^*)(1-\alpha^*\tau)+
\beta_{\lambda}(\lambda^*,\mu^*)(\tau\beta^{*2}-\alpha^*)}&=&0,\label{cond2}
\end{eqnarray}
then the path $\lambda\longmapsto (\alpha(\lambda,\mu),\beta(\lambda,\mu))$ has a tangential intersection with the Hopf bifurcation curve (\ref{bifcurve}) at the point $(\alpha^*,\beta^*)$.
At the point $(\alpha^*,\beta^*)$, the signed curvature of the Hopf bifurcation curve (\ref{bifcurve}) can be computed as
\[
\kappa_1=\frac{\beta^{*}(\alpha^{*2}-\beta^{*2})(\beta^{*2}\tau^2+\alpha^*\tau-2)\tau}{[(\beta^{*2}\tau^2+1)(\alpha^{*2}+\beta^{*2})-4\alpha^*\beta^{*2}\tau]^{3/2}}
\]
and the signed curvature of the path $\lambda\longmapsto (\alpha(\lambda,\mu),\beta(\lambda,\mu))$ is
\[
\kappa_2=\frac{(\alpha^*\tau-1)^2\beta^{*2}(\beta^*(\alpha^*\tau-1)\alpha_{\lambda\lambda}(\lambda^*,\mu^*)+(\alpha^*-\beta^{*2}\tau)\beta_{\lambda\lambda}(\lambda^*,\mu^*))}{\beta_{\lambda}(\lambda^*,\mu^*)^2[(\beta^{*2}\tau^2+1)(\alpha^{*2}+\beta^{*2})-4\alpha^*\beta^{*2}\tau]^{3/2}}.
\]
Although these two curves are tangential at the intersection point, we want their curvatures to be different, so that locally near $(\alpha^*,\beta^*)$ the curve  $\lambda\longmapsto (\alpha(\lambda,\mu),\beta(\lambda,\mu))$ lies entirely on one side only of the Hopf curve, as illustrated in Figure \ref{fig3}.  This is equivalent to requiring
\begin{equation}
\begin{array}{r}
{\cal G}(\alpha^*,\beta^*,\tau)\equiv (\alpha^*\tau-1)^3\beta^{*2}\alpha_{\lambda\lambda}(\lambda^*,\mu^*)+\beta^*(\alpha^*\tau-1)^2(\alpha^*-\beta^{*2}\tau)\beta_{\lambda\lambda}(\lambda^*,\mu^*)\\[0.2in]
-(\beta_{\lambda}(\lambda^*,\mu^*))^2\tau (\alpha^{*2}-\beta^{*2})(\beta^{*2}\tau^2+\alpha^*\tau-2)\neq 0.
\end{array}
\label{cond3}
\end{equation}

By using center manifold and normal form analysis for retarded functional differential equations (see for example \cite{FM1,FM2}), we expect to see for $(\lambda,\mu)$ near $(\lambda^*,\mu^*)$ Hopf bifurcation diagrams for 
\[
\dot{x}(t)=\alpha(\lambda,\mu)\,x(t)+\beta(\lambda,\mu)\,x(t-\tau)+F(x(t),x(t-\tau),\lambda,\mu)
\] 
which resemble those of Figure \ref{fig5}.  This will be the theoretical focus of section 2 of this paper.

\subsection{Endemic bubbles}

In \cite{bubbles}, the following susceptible-infected-susceptible (SIS) model incorporating delayed behavioral response was analyzed
\begin{equation}
\dot{y}(t)=-y(t)+R_0h(y(t-\tau))y(t)(1-y(t))
\label{SIS}
\end{equation}
where $y(t)$ represents the proportion of infected individuals in the population, $R_0$ is the basic reproduction number (expressing the expected number of secondary infections generated by a single infectious agent introduced into a wholly susceptible population), and the smooth behavioral function $h:[0,1]\longrightarrow (0,1]$ is such that $h'(y)\leq 0$, $h(0)=1$ and $h(1)<1$.  The parameter $\tau>0$ represents the delay in time 
between the moment when the population has information about the number of infections, and the moment this population decides to modify its behavior in order to reduce the rate of infections.
We refer the reader to \cite{bubbles} for more details related to the modelling aspects.

\begin{figure}[htpb]
\begin{center}
\includegraphics[width=3in]{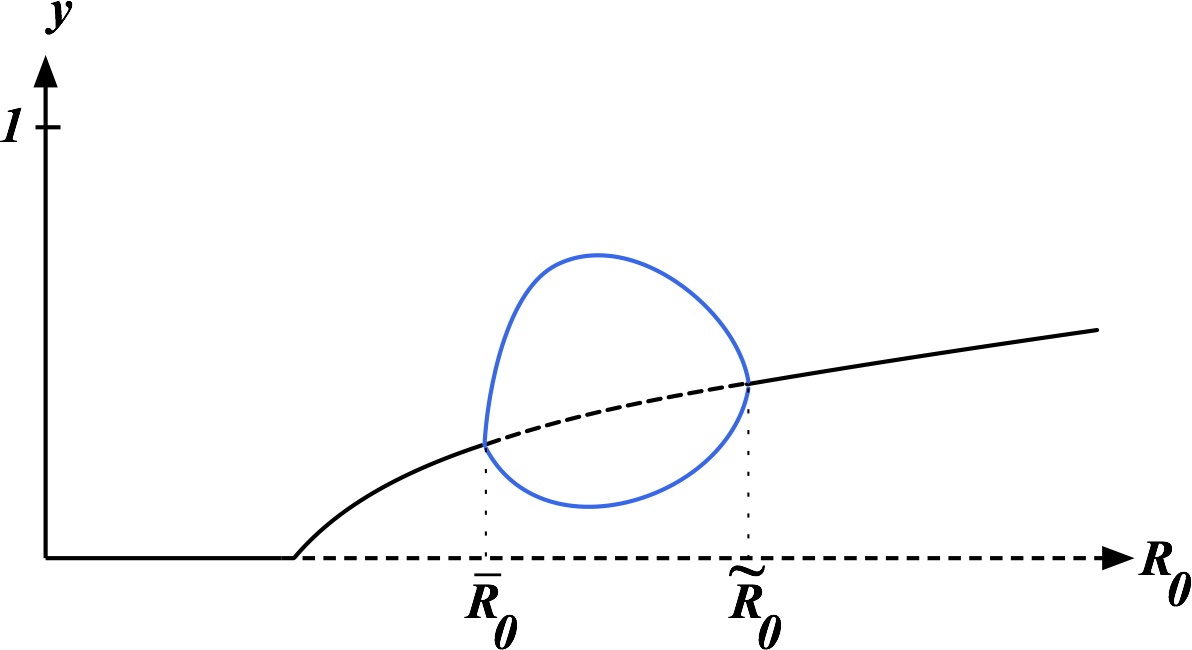}
\caption{Schematic representation of a bifurcation diagram for (\ref{SIS}) which exhibits the phenomenon of {\em endemic bubble}.  The disease-free equilibrium is stable for $R_0\leq 1$.  At $R_0=1$, a branch of endemic equilibria bifurcate from the disease-free equilibrium, and this branch is stable until $R_0=\overline{R}_0$ at which point a Hopf bifurcation occurs.  The blue curves represent minima and maxima of oscillations for the bifurcating periodic solutions.  At $R_0=\tilde{R}_0$, the branch of periodic solutions reconnects with the endemic equilibrium.  The region between $\overline{R}_0$ and $\tilde{R}_0$ is what is referred to as the {\em endemic bubble}.}
\label{fig6}
\end{center}
\end{figure}

Equilibria of (\ref{SIS}) are algebraic solutions to
\begin{equation}
y=R_0h(y)y(1-y).
\label{eqeq}
\end{equation}
Obviously $y=0$ is always a solution to (\ref{eqeq}) (disease-free equilibrium).  It is shown in \cite{bubbles} that if $R_0\leq 1$, then the disease-free equilibrium is globally asymptotically stable.  If $R_0>1$, then $y=0$ becomes unstable, and (\ref{eqeq}) has a unique {\em endemic equilibrium} $y^*$ satisfying $y^*<1-1/R_0$.  The authors in \cite{bubbles} then perform a comprehensive analysis of (\ref{SIS}) for different choices of response functions $h$, and in particular, they plot bifurcation diagrams with distinguished parameter $R_0$.  Some of these bifurcation diagrams (see for example figures 5 and 9 of \cite{bubbles}) exhibit a phenomenon that the authors have called {\em endemic bubbles}, which loosely speaking is the bifurcation of a branch of periodic solutions from the endemic equilibrium point at some parameter value $\overline{R}_0>1$, and this branch reconnects with the endemic equilibrium (in a reverse Hopf bifurcation) at parameter value $\tilde{R}_0>\overline{R}_0$.  A schematic representation of a typical such endemic bubble bifurcation diagram is given in Figure \ref{fig6}.  

When focussing on the endemic branch, we notice the similarity between Figure \ref{fig6} and the bottom-right-most bifurcation in Figure \ref{fig5}.  As we will show in this paper, the endemic bubble is indeed a consequence of a degenerate Hopf bifurcation which occurs in (\ref{SIS}) for various choices of response functions $h(y,p)$, where $p\geq 0$ is some auxiliary parameter which appears in the behavioral response function. 

\subsection{Overview of paper}

With the above discussion in mind, this suggests that the organizing center of endemic bubbles is a degenerate Hopf bifurcation in (\ref{SIS}).  In this paper, we will develop the theoretical ideas to make this hypothesis into a rigorous theorem.  Our approach will be to use the center manifold and normal form theory for RFDEs developed in \cite{FM1,FM2} and the classification and unfolding of degenerate Hopf bifurcations of \cite{GoLang}. 

Ours is not the first study of higher codimension bifurcation in RFDEs with nonlinear degeneracies.  We note in particular \cite{BR} where the authors present and analyze a DDE model for tumor growth in which a Bautin bifurcation (Hopf bifurcation with the first Lyapunov coefficient being zero) occurs.  
In \cite{Qesmi}, a Maple program is presented which allows the computation of coefficients up to any order related to Hopf bifurcation for RFDEs.
Also, in \cite{BPS}, the authors investigate the stabilization of unstable periodic orbits via the Pyragas delayed feedback control.  It is shown that the stabilization mechanism occurs because of a highly degenerate Hopf bifurcation in which both crossing condition and first Lyapunov coefficient degeneracy occur simultaneously. 
Our paper here appears to be the first systematic study of 
the degeneracy resulting from the violation of the crossing condition in the Hopf Bifurcation Theorem for general RFDEs.

In section 2, we will establish the sufficient conditions for such a degenerate Hopf bifurcation to occur in a general class of parametrized RFDEs.  Our main theorems to that effect are Theorems \ref{Mainthm23} and \ref{thmmain}.  In subsection 2.3, we will apply these results to the special case in which the RFDE is a DDE of the form (\ref{lindde1}).

In section 3, we use our theoretical results to seek out and analyze degenerate Hopf bifurcation points in the SIS model (\ref{SIS}), where the behavioral response function $h=h(y,p)$ depends on an auxiliary parameter $p$.  The parameter space is thus two-dimensional: $R_0$ acts as a distinguished bifurcation parameter, and $p$ as an unfolding parameter.  This leads to Theorem \ref{thm42}.  In subsections 3.2 and 3.3, we illustrate these results by performing explicit computations for two of the types of delayed response functions $h(y,p)$ which were considered in \cite{bubbles}.  

We end with some concluding remarks in section 4, and an Appendix where a lengthy expression for the Lyapunov coefficient associated with the degenerate Hopf bifurcation is presented.

\Section{Center manifold and normal form reduction}

In this section, we give a brief summary of the theory presented in \cite{FM1,FM2} for center-manifold and normal form reductions of retarded functional differential equations, and then apply it to study the degenerate Hopf bifurcation in a class of parametrized RFDEs.

\subsection{Phase space and splitting}

For $\tau>0$, we consider the Banach space $C=C([-\tau,0],\mathbb{R})$ of continuous functions from $[-\tau,0]$ into $\mathbb{R}$, endowed with the supremum norm.  We define $z_t\in C$ by $z_t(\theta)=z(t+\theta)$, for $\theta\in [-\tau,0]$.

Let 
\begin{equation}
\dot{z}(t)=L(\lambda,\mu)(z_t)+F(z_t,\lambda,\mu)
\label{RFDE1}
\end{equation}
denote a smoothly parametrized family of nonlinear retarded functional differential equations, where $\lambda\in\mathbb{R}$ is a real distinguished bifurcation parameter, $\mu\in\mathbb{R}$ is an auxiliary parameter (to be regarded as an unfolding parameter), $L(\lambda,\mu)$ is a smoothly parametrized family of bounded linear operators from $C$ into $\mathbb{R}$, and $F$ is a smooth function from $C\times\mathbb{R}^2$ into $\mathbb{R}$, which satisfies
\begin{equation}
F(0,\lambda,\mu)=0,\,\,\,\,D_1F(0,\lambda,\mu)=0,\,\,\,\,\forall (\lambda,\mu)\,\,\,\mbox{\rm near}\,\,\,(0,0)\,\,\in\mathbb{R}^2
\label{zerotayl}
\end{equation}
We denote $L_0=L(0,0)$, and we rewrite (\ref{RFDE1}) as
\[
\dot{z}(t)=L_0(z_t)+(L(\lambda,\mu)-L_0)(z_t)+F(z_t,\lambda,\mu).
\]

By the Riesz representation theorem, we may write
\[
L(\lambda,\mu)(\phi)=\int_{-\tau}^0 [d\eta_{\lambda,\mu}(\theta)]\,\phi(\theta),
\]
where $\eta_{\lambda,\mu}(\theta)$ is a measurable function on $[-\tau,0]$.  We also define ${\cal A}(\lambda,\mu)$ to be the infinitesimal generator for the flow of the linear system $\dot{z}=L(\lambda,\mu)(z_t)$.
For each fixed $(\lambda,\mu)$, we consider the characteristic equation
\[
\Delta(\lambda,\mu)(\xi)=0,\,\,\,\,\,\,\Delta(\lambda,\mu)(\xi)=\xi-\int_{-\tau}^0\,d\eta_{\lambda,\mu}(\theta)\,e^{\xi\theta}.
\]
We suppose that the following holds:
\begin{hyp} The characteristic equation $\Delta(\lambda,\mu)(\xi)=0$ has roots, $\xi(\lambda,\mu)=\gamma(\lambda,\mu)\pm i\omega(\lambda,\mu)$ which are smoothly parametrized by $\lambda$ and $\mu$, and are such that
\[
\gamma(0,0)=0,\,\,\,\,\omega(0,0)\equiv\omega^* > 0,\,\,\,\,\,\gamma_{\lambda}(0,0)=0,\,\,\,\,\,\,\gamma_{\mu}(0,0)\neq 0,
\]
\begin{equation}
\frac{d}{d\xi}\Delta(0,0)(i\omega^*)\neq 0,\,\,\,\,\mbox{\rm or equivalently}\,\,\,\,\,1-L_0(\theta\,e^{i\omega^*\theta})\neq 0,
\label{simplecondition}
\end{equation}
and furthermore, we suppose the characteristic equation $\Delta(0,0)(\xi)=0$ has no roots other than $\pm\,i\omega^*$ on the imaginary axis, and that all other roots of $\Delta(0,0)(\xi)=0$ have strictly negative real part.
\label{hyp1}
\end{hyp}

It follows that if ${\cal A}_0$ denotes the infinitesimal generator of the linear equation $\dot{z}(t)=L_0(z_t)$, then the eigenspace $P$ of ${\cal A}_0$ corresponding to the eigenvalues $\pm i\omega^*$ is two-dimensional.
Let $\Phi(\theta)=(e^{i\omega^*\theta},e^{-i\omega^*\theta})$ be a complex representation of a basis for $P$.    The dual space $C^*=C([0,\tau],\mathbb{R})$ also admits a two-dimensional subspace for the transpose ${\cal A}_0^*$ corresponding to $\pm i\omega^*$.  We introduce the bilinear form between $C$ and $C^*$
\[
(\psi,\phi)=\psi(0)\phi(0)-\int_{-\tau}^0\int_0^\theta\,\psi(\upsilon-\theta)\,d\eta(\theta)\,\phi(\upsilon)\,d\upsilon.
\]
We choose a basis $\Psi(s)=\mbox{\rm col}(\psi_1(0)e^{-i\omega s},\overline{\psi_1(0)}e^{i\omega s})$ such that $(\Psi,\Phi)=I_d$.  As is shown in \cite{FM2} we have
\begin{equation}
\psi_1(0)=[1-L_0(\theta e^{i\omega^*\theta})]^{-1},
\label{psidef}
\end{equation}
which is well-defined because of (\ref{simplecondition}).
We have $C=P\oplus Q$, where $Q$ is infinite-dimensional and also invariant for the operator ${\cal A}_0$.

We now let $BC$ represent the Banach space of functions from $[-\tau,0]$ into $\mathbb{R}$ which are uniformly continuous on $[-\tau,0)$ with a jump discontinuity at 0.  We define the function
\[
X_0(\theta)=\left\{\begin{array}{ll}1,&\theta=0\\[0.1in]
0,&-\tau\leq\theta<0.
\end{array}
\right.
\]
The elements of the space $BC$ can be written as $Y=\varphi+X_0\nu$, where $\varphi\in C$ and $\nu\in\mathbb{R}$.  We define the projection $\pi: BC\longrightarrow P$ as
\[
\pi(\varphi+X_0\nu)=\Phi [ (\Psi,\varphi)+\Psi(0)\nu].
\]
We then have
\[
BC=P\oplus \mbox{\rm ker}\,\pi,
\]
where $Q\subsetneq \mbox{\rm ker}\,\pi$, and we may rewrite (\ref{RFDE1}) according to this splitting as
\begin{equation}
\begin{array}{c}
\dot{x}=Bx+\Psi(0)[(L(\lambda,\mu)-L_0)(\Phi\,x+\zeta)+F(\Phi\,x+\zeta,\lambda,\mu)],\\[0.1in]
\frac{d}{dt}\zeta=A_{Q_1}\zeta+(I_d-\pi)X_0[(L(\lambda,\mu)-L_0)(\Phi\,x+\zeta)+F(\Phi\,x+\zeta,\lambda,\mu)],
\end{array}
\label{split}
\end{equation}
where 
\[
B=\left(\begin{array}{cc}i\omega^*&0\\0&-i\omega^*\end{array}\right),
\]
$x=(u,\overline{u})^T$, $\zeta\in Q^1\equiv Q\cap C^1([-\tau,0],\mathbb{R})$, and $A_{Q_1}$ is defined as
\[
A_{Q_1}\varphi=\dot{\varphi}+X_0 [ L_0\varphi-\dot{\varphi}(0)].
\]

As is shown in \cite{FM2}, it is then possible to define near identity changes of coordinates
\begin{equation}
(x,\zeta)=(\hat{x},\hat{\zeta})+U(\hat{x},\lambda,\mu)
\label{coc}
\end{equation}
such that a Taylor expansion of (\ref{split}) can be put into the normal form
\[
\dot{x}=Bx+\sum_{j\geq 2}\,g_j^1(x,\zeta,\lambda,\mu),\,\,\,\,\,\,\,\frac{d}{dt}\zeta=A_{Q_1}\zeta+\sum_{j\geq 2}\,g_j^2(x,\zeta,\lambda,\mu)
\]
having the property that the center manifold is locally given by $\zeta=0$, and the local flow of (\ref{RFDE1}) on this center manifold is given by
\begin{equation}
\dot{x}=Bx+\sum_{j\geq 2}\,g_j^1(x,0,\lambda,\mu).
\label{cmnf}
\end{equation}
The nonlinear terms in (\ref{cmnf}) are in classical Poincar\'e-Birkhoff normal form with respect to the matrix $B$.

\subsection{Normal form for the degenerate Hopf bifurcation}

From Hypothesis \ref{hyp1}, we have that there exists a smooth function $\xi(\lambda,\mu)$ such that
\[
\xi(\lambda,\mu)=L(\lambda,\mu)(e^{\xi(\lambda,\mu)\theta}),\,\,\,\,\,\,\xi(0,0)=i\omega^*,
\]
for which implicit differentiation gives the following equalities
\[
\xi_{\lambda}(0,0)=\psi_1(0)\,L_{\lambda}(0,0)(e^{i\omega^*\theta}),
\]
\[
\xi_{\mu}(0,0)=\psi_1(0)\,L_{\mu}(0,0)(e^{i\omega^*\theta}),
\]
and
\begin{equation}
\xi_{\lambda\lambda}(0,0)=\psi_1(0)\,\left[L_{\lambda\lambda}(0,0)(e^{i\omega^*\theta})+2\xi_{\lambda}(0,0)L_{\lambda}(0,0)(\theta\,e^{i\omega^*\theta})+(\xi_{\lambda}(0,0))^2L_0(\theta^2\,e^{i\omega^*\theta})\right],
\label{xilambdalambda}
\end{equation}
where $\psi_1(0)$ is as in (\ref{psidef}).
Hypothesis \ref{hyp1} thus implies
\[
\mbox{\rm Re}\left[ \psi_1(0)\,L_{\lambda}(0,0)(e^{i\omega^*\theta}) \right] = 0,\,\,\,\,\,\,\,\mbox{\rm and}\,\,\,\,\,\,
\mbox{\rm Re}\left[ \psi_1(0)\,L_{\mu}(0,0)(e^{i\omega^*\theta}) \right]\neq 0.
\]
In addition, we are going to assume the following second-order non-degeneracy condition
\begin{hyp}
We assume the non-degeneracy condition $\mbox{\rm Re} (\xi_{\lambda\lambda}(0,0))\neq 0$, which is equivalent (via (\ref{xilambdalambda}))
to
\begin{equation}
\mbox{\rm Re}\left(\psi_1(0)\,\left[L_{\lambda\lambda}(0,0)(e^{i\omega^*\theta})+2\xi_{\lambda}(0,0)L_{\lambda}(0,0)(\theta\,e^{i\omega^*\theta})+(\xi_{\lambda}(0,0))^2L_0(\theta^2\,e^{i\omega^*\theta})\right]\right)\neq 0.
\label{nondeg2}
\end{equation}
\label{hyp2}
\end{hyp}
We will see later on that the geometrical meaning of (\ref{nondeg2}) is equivalent to the curvature condition given in (\ref{cond3}).

Consider equations (\ref{split}) which we expand in a Taylor series as
\begin{equation}
\begin{array}{rcl}
\dot{x}&=&Bx+\Psi(0)[\lambda\,L_{\lambda}(0,0)(\Phi\,x+\zeta)+\mu\,L_{\mu}(0,0)(\Phi\,x+\zeta)\\[0.1in]
&&+\frac{1}{2}\lambda^2\,L_{\lambda\lambda}(0,0)(\Phi\,x+\zeta)+R(x,\zeta,\lambda,\mu)],\\[0.2in]
\frac{d}{dt}\zeta&=&A_{Q_1}\zeta+(I_d-\pi)X_0[\lambda\,L_{\lambda}(0,0)(\Phi\,x+\zeta)+\mu\,L_{\mu}(0,0)(\Phi\,x+\zeta)\\[0.1in]
&&+\frac{1}{2}\lambda^2\,L_{\lambda\lambda}(0,0)(\Phi\,x+\zeta)+R(x,\zeta,\lambda,\mu)],
\end{array}
\label{expan}
\end{equation}
where the remainder term is 
$R(x,\zeta,\lambda,\mu)=O(|(x,\zeta)|^2,\mu^2,\mu\lambda,\lambda^3)$.

We procede as in section 3 of \cite{FM2}. Writing $x=(u,\overline{u})^T$, $\Phi\,x=ue^{i\omega^*\theta}+\bar{u}e^{-i\omega^*\theta}$, we define quadratic near-identity changes of coordinates of the form (\ref{coc}) which transforms (\ref{expan}) in such a way that the normal form (\ref{cmnf}) becomes
\begin{equation}
\dot{u}=i\omega^*u+[(\sigma_1+i\sigma_2)\mu+i\sigma_3\lambda+(\sigma_4+i\sigma_5)\lambda^2]u+\psi_1(0)R((u,\overline{u}),0,\lambda,\mu),
\label{TE}
\end{equation}
where
\[
\sigma_1+i\sigma_2=\psi_1(0)\,L_{\mu}(0,0)(e^{i\omega^*\theta}),\,\,\,\,\,\sigma_1\neq 0,
\]
\[
\sigma_3=\mbox{\rm Im}\left[\psi_1(0)L_{\lambda}(0,0)(e^{i\omega^*\theta})\right]=\omega_{\lambda}(0,0),
\]
and
\begin{equation}
\begin{array}{rcl}
\sigma_4+i\sigma_5&=&\frac{1}{2}\psi_1(0)\,\left[L_{\lambda\lambda}(0,0)(e^{i\omega^*\theta})+2\xi_{\lambda}(0,0)L_{\lambda}(0,0)(\theta\,e^{i\omega^*\theta})+(\xi_{\lambda}(0,0))^2L_0(\theta^2\,e^{i\omega^*\theta})\right]\\[0.2in]
&=&\frac{1}{2}\xi_{\lambda\lambda}(0,0)
\end{array}
\label{sig45def}
\end{equation}
is such that $\sigma_4\neq 0$ because of Hypothesis \ref{hyp2}.  The extra terms 
\[
\frac{1}{2}\psi_1(0)\,\left[2\xi_{\lambda}(0,0)L_{\lambda}(0,0)(\theta\,e^{i\omega^*\theta})+(\xi_{\lambda}(0,0))^2L_0(\theta^2\,e^{i\omega^*\theta})\right]
\]
which appear in (\ref{sig45def}) arise from having normalized quadratic terms in (\ref{expan}).

As mentioned earlier, the nonlinear terms in (\ref{cmnf}) (equivalently (\ref{TE})) are in classical Poincar\'e-Birkhoff normal form with respect to the matrix $B$.  It is well-known (see for example \cite{GSS88}) that such a normal form has the algebraic form
\begin{equation}
\dot{u}=i\omega^*u+\Gamma(|u|^2,\lambda,\mu)u,
\label{nf22}
\end{equation}
where $\Gamma$ is a smooth complex-valued function.

The following theorem now follows from the above discussion and by performing further near-identity changes of coordinates (\ref{coc}) as in \cite{FM2}:
\begin{thm}
Consider the smoothly parametrized family of nonlinear retarded functional differential equations (\ref{RFDE1}) which satisfies (\ref{zerotayl}), and Hypotheses \ref{hyp1} and \ref{hyp2}.  Then there exists a two-dimensional semi-flow invariant center manifold in the phase space $C([-\tau,0],\mathbb{R})$.  Furthermore, there exist a formal sequence of parameter-dependent near-identity changes of coordinates of the form (\ref{coc}) which are such that the dynamics of (\ref{RFDE1}) reduced to this center manifold have a Poincar\'e-Birkhoff normal form to any order given by (\ref{cmnf}) (equivalently (\ref{nf22})):
\begin{equation}
\dot{u}=\left[i\omega^*+(\sigma_1+i\sigma_2)\mu+i\sigma_3\lambda+(\sigma_4+i\sigma_5)\lambda^2+H(|u|^2,\lambda,\mu)\right] u,
\label{nf44}
\end{equation}
where $H(|u|^2,0,0)\equiv K_1+iK_2$, $H_\lambda(|u|^2,0,0)=H_{\mu}(|u|^2,0,0)=H_{\lambda\lambda}(|u|^2,0,0)=0$.
\label{Mainthm23}
\end{thm}
The coefficient $K_1$ in (\ref{nf44}) is the first Lyapunov coefficient.  It has been computed explicitly in \cite{FM2} as
\begin{equation}
K_1=\mbox{\rm Re}\left[\psi_1(0)\left(B_{(2,1,0,0)}-\frac{B_{(1,1,0,0)}B_{(1,0,1,0)}}{L_0(1)}+\frac{B_{(2,0,0,0)}B_{(0,1,0,1)}}{2i\omega^*-L_0(e^{2i\omega^*\theta})}\right)\right],
\label{Kdef}
\end{equation}
where the coefficients $B_{(i,j,k,\ell)}$ are read-off from the Taylor expansion of $F$ in (\ref{RFDE1}):
\begin{equation}
\begin{array}{l}
F(x_1e^{i\omega^*\theta}+x_2e^{-i\omega^*\theta}+x_3\,1+x_4e^{2i\omega^*\theta},0,0)=\\[0.2in]
B_{(2,0,0,0)}x_1^2+B_{(1,1,0,0)}x_1x_2+B_{(1,0,1,0)}x_1x_3+B_{(0,1,0,1)}x_2x_4+B_{(2,1,0,0)}x_1^2x_2+\cdots
\end{array}
\label{Kdef2}
\end{equation}
Using the polar coordinates $u=re^{i\phi}$, (\ref{nf44}) becomes
\begin{equation}
\begin{array}{rcl}
\dot{r}&=&r(\sigma_4\lambda^2+\sigma_1\mu+K_1r^2+O(r^4,\mu^2,\mu\lambda,\lambda^3))\\[0.2in]
\dot{\phi}&=&\omega^*+\cdots
\end{array}
\label{polarnf}
\end{equation}
Zeros of the $\dot{r}$ equation of (\ref{polarnf}) correspond to periodic solutions of (\ref{RFDE1}) whose period is approximately equal to $2\pi/\omega^*$.  Let us assume the non-degeneracy condition $K_1\neq 0$.
Since $\sigma_1\neq 0$ and $\sigma_4\neq 0$, it now follows from the classification of degenerate Hopf bifurcations done in \cite{GoLang} that
\begin{thm}
For $\mu$ near $0$, the bifurcation diagram (wrt $\lambda$) of 
$r(\sigma_4\lambda^2+\sigma_1\mu+K_1r^2+O(r^4,\mu^2,\mu\lambda,\lambda^3))=0$ is locally diffeomorphic to the bifurcation diagram of
\[
r(\varepsilon(\lambda^2+\eta)+r^2)=0,
\]
where $\varepsilon=\mbox{\rm sgn}(\sigma_4/K_1)$ and 
${\displaystyle\eta=\frac{\sigma_1\mu}{|K_1|\mbox{\rm sgn}(\sigma_4)}}$ (see (\ref{normform})), as illustrated in Figure \ref{fig5}. 
\label{thmmain}
\end{thm}

\subsection{Special case: discrete delay}

As a special case of the theory we have just developed, we return to the prototype delay-differential equation (\ref{lindde1})
\begin{equation}
\dot{x}(t)=\alpha(\lambda,\mu)x(t)+\beta(\lambda,\mu)x(t-\tau)+F(x(t),x(t-\tau),\lambda,\mu)
\label{prot}
\end{equation}
where we will assume without loss of generality that $\lambda^*=0$, $\mu^*=0$, and $(\alpha^*,\beta^*)$ are such that (\ref{cond1}), (\ref{cond2}) and (\ref{cond3}) are satisfied, and $\beta^*\neq 0$, and $\alpha^*\tau-1\neq 0$.

The operator $L(\lambda,\mu)$ is defined as
\[
L(\lambda,\mu)(z(\theta))=\alpha(\lambda,\mu)z(0)+\beta(\lambda,\mu)z(-\tau),
\]
so that
\[
\begin{array}{rcl}
L_0(e^{i\omega^*\theta})&=&{\displaystyle\alpha^*+\beta^*e^{-i\omega^*\tau}=i\omega^*\,\,\Longrightarrow e^{-i\omega^*\tau}=\frac{i\omega^*-\alpha^*}{\beta^*}},\\[0.2in]
L_0(\theta\,e^{i\omega^*\theta})&=&-\tau \beta^* e^{-i\omega^*\tau}=-{\tau}(i \omega^*-\alpha^*),\\[0.2in]
\psi_1(0)&=&[1-L_0(\theta\,e^{i\omega^*\theta})]^{-1}={\displaystyle\frac{1}{(1-\alpha^*\tau)+i\omega^*\tau}}.
\end{array}
\]
Using (\ref{cond2}), we have
\[
\begin{array}{ll}
\psi_1(0)L_{\lambda}(0,0)(e^{i\omega^*\theta})&{\displaystyle=\frac{\alpha_{\lambda}(0,0)+\frac{\beta_{\lambda}(0,0)}{\beta^*}(i\omega^*-\alpha^*)}{(1-\alpha^*\tau)+i\omega^*\tau}}\\[0.2in]
&{\displaystyle=-i\left[\frac{\omega^*(\alpha_{\lambda}(0,0)\beta^*\tau-\beta_{\lambda}(0,0))}{\beta^*((1-\alpha^*\tau)^2+\omega^{*2}\tau^2)}\right].}
\end{array}
\]
The non-degeneracy condition $\mbox{\rm Re}(\psi_1(0)L_{\mu}(0,0)(e^{i\omega^*\theta}))$ becomes
\[
\sigma_1\equiv \frac{{\displaystyle\beta^*\alpha_{\mu}(0,0)(1-\alpha^*\tau)+
\beta_{\mu}(0,0)(\tau\beta^{*2}-\alpha^*)}}{\beta^*((1-\alpha^*\tau)^2+\omega^{*2}\tau^2)}\neq 0.
\]
Finally,
\[
\begin{array}{rcl}
2\sigma_4&=&\mbox{\rm Re}(\psi_1(0)[L_{\lambda\lambda}(0,0)(e^{i\omega^*\theta})+2\xi_{\lambda}(0,0)L_{\lambda}(0,0)(\theta\,e^{i\omega^*\theta})+(\xi_{\lambda}(0,0))^2L_0(\theta^2\,e^{i\omega^*\theta})])\\[0.1in]
&=&{\displaystyle\frac{{\cal G}(\alpha^*,\beta^*,\tau)}{\beta^{*2}(\alpha^*\tau-1)^2(-\beta^{*2}\tau^2+2\alpha^*\tau-1)}},
\end{array}
\]
%
%
where ${\cal G}(\alpha^*,\beta^*,\tau)$ is as in (\ref{cond3}).  Thus, we note that the condition $\sigma_4\neq 0$ is equivalent to the curvature condition (\ref{cond3}).

We write $F$ in (\ref{prot}) at $(\lambda,\mu)=(0,0)$ as
\begin{equation}
\begin{array}{ll}
F(x(t),x(t-\tau),0,0)&=f_{(2,0)}x(t)^2+f_{(1,1)}x(t)x(t-\tau)+f_{(0,2)}(x(t-\tau))^2+\\[0.2in]
&+f_{(3,0)}x(t)^3+f_{(2,1)}x(t)^2x(t-\tau)+f_{(1,2)}x(t)(x(t-\tau))^2+f_{(0,3)}(x(t-\tau))^3\\[0.2in]
&+O(|x|^4)
\end{array}
\label{Fdef}
\end{equation}
To compute the first Lyapunov coefficient $K_1$ in (\ref{Kdef}), we need
\[
L_0(1)=\alpha^*+\beta^*,\,\,\,\,\,\mbox{\rm and}\,\,\,\,\,L_0(e^{2i\omega^*\theta})=\alpha^*+\beta^*e^{-2i\omega^*\tau}=\alpha^*+\frac{(i\omega^*-\alpha^*)^2}{\beta^*}.
\]
We will assume that $\alpha^*+\beta^*\neq 0$.  If this condition holds, it follows that $2i\omega^*-L_0(e^{2i\omega^*\theta})\neq 0$.
A lengthy computation using (\ref{Kdef}), (\ref{Kdef2}) and (\ref{Fdef}) yields the Lyapunov coefficient $K_1$ in terms of the coefficients $f_{(j,k)}$ in (\ref{Fdef}) and $\alpha^*$, $\beta^*$, $\omega^*$ and $\tau$.  The formula is lengthy and given in the Appendix.  Generically, we will have $K_1\neq 0$.  Once we have these quantities, we can compute the
normal form for the unfolding of the degenerate Hopf bifurcation at $(\lambda,\mu)=(0,0)$:
\[
r(\varepsilon (\lambda^2+\eta)+r^2)=0,
\]
where $\varepsilon=\mbox{\rm sgn}(\sigma_4/K_1)$ and ${\displaystyle\eta=\frac{\sigma_1\mu}{|K_1|\mbox{\rm sgn}(\sigma_4)}}$ as in Theorem \ref{thmmain}.

\Section{Application to a SIS model with delayed behavioral response: endemic bubbles}

We recall the SIS model (\ref{SIS}) which was analyzed in \cite{bubbles}.  In this section, we will apply the theoretical results of the previous section to establish and study degenerate Hopf bifurcations which occur in this model.

\subsection{General case}

Non-trivial equilibria of (\ref{SIS}) are solutions $\bar{y}\neq 0$ to (\ref{eqeq}), or equivalently to
\begin{equation}
h(\bar{y},p)=\frac{1}{R_0 (1-\bar{y})}.
\label{eqeq2}
\end{equation}
Because the function $h$ is such that $h_y(y,p)\leq 0$, $h(0,p)=1$ and $h(1,p)<1$, if $R_0\geq 1$ then (\ref{eqeq2}) has a unique solution $\bar{y}=\bar{y}(R_0,p)$, i.e.
\begin{equation}
R_0h(\bar{y}(R_0,p),p)(1-\bar{y}(R_0,p))\equiv 1,\,\,\,\,\,\,\,\,\,\,\,\forall R_0\geq 1, p\geq 0.
\label{eqeq3}
\end{equation}
Linearizing (\ref{SIS}) about the equilibrium $\bar{y}(R_0,p)$ gives
\[
\dot{x}(t)=\alpha(R_0,p)x(t)+\beta(R_0,p)x(t-\tau),
\]
where
\begin{equation}
\alpha(R_0,p)=-R_0\,\bar{y}(R_0,p)\,h(\bar{y}(R_0,p),p)<0
\label{linalpha}
\end{equation}
and
\begin{equation}
\beta(R_0,p)=R_0\,h_y(\bar{y}(R_0,p),p)\,\bar{y}(R_0,p)\,(1-\bar{y}(R_0,p))<0.
\label{linbeta}
\end{equation}
Using implicit differentiation of (\ref{eqeq3}), one can compute $\bar{y}_{R_0}$, $\bar{y}_{p}$ and $\bar{y}_{R_0R_0}$ in terms of $h$, $h_y$, $h_{yy}$, $h_p$ and $\bar{y}$.  One can then use (\ref{linalpha}) and (\ref{linbeta}) to compute the quantities $\alpha_{R_0}$, $\beta_{R_0}$, $\alpha_p$, $\beta_p$, $\alpha_{R_0R_0}$ and $\beta_{R_0R_0}$.

We then have the following
\begin{thm}
Consider the nonlinear delay-differential equation
\[
\dot{x}(t)=\alpha(R_0,p)x(t)+\beta(R_0,p)x(t-\tau)+{\cal F}(x(t),x(t-\tau),R_0,p)
\]
which is obtained by performing the change of variables $x=y-\bar{y}(R_0,p)$ in (\ref{SIS}), (where $\bar{y}(R_0,p)$ is the endemic equilibrium, which is solution to (\ref{eqeq3})), and where $\alpha(R_0,p)$ and $\beta(R_0,p)$ are as in (\ref{linalpha}) and (\ref{linbeta}) respectively.  Define $\omega(R_0,p)$ by
\[
\omega(R_0,p)^2=\beta(R_0,p)^2-\alpha(R_0,p)^2.
\]
Suppose $(R_0^*,p^*)$ is a point in parameter space such that
\begin{eqnarray}
\alpha^*+\beta^*\cos\omega^*\tau&=&0\label{condeq1},\\[0,2in]
\beta^*\alpha_{R_0}^*\,(1-\alpha^*\tau)+\beta_{R_0}^*\,(\tau\beta^{*2}-\alpha^*)&=&0,\label{condeq2}
\end{eqnarray}
where $\alpha^*\equiv \alpha(R_0^*,p^*)$, $\beta^*\equiv\beta(R_0^*,p^*)$, $\alpha_{R_0}^*\equiv \alpha_{R_0}(R_0^*,p^*)$, $\beta_{R_0}^*\equiv\beta(R_0^*,p^*)$, $\alpha_p^*\equiv\alpha_p(R_0^*,p^*)$, $\beta_p^*\equiv\beta_p(R_0^*,p^*)$, $\alpha_{R_0R_0}^*\equiv\alpha_{R_0R_0}(R_0^*,p^*)$,
$\beta_{R_0R_0}^*\equiv\beta_{R_0R_0}(R_0^*,p^*)$, and $\omega^{*2}\equiv \beta^{*2}-\alpha^{*2}>0$.  Suppose
\begin{equation}
\sigma_1^*\equiv\frac{\beta^*\alpha_{p}^*\,(1-\alpha^*\tau)+\beta_{p}^*\,(\tau\beta^{*2}-\alpha^*)}{\beta^*((1-\alpha^*\tau)^2+\omega^{*2}\tau^2)}\neq 0,
\label{condeq10}
\end{equation}
\begin{equation}
2\sigma_4^*\equiv\frac{{\cal G}(\alpha^*,\beta^*,\tau)}{\beta^{*2}(\alpha^*\tau-1)^2(-\beta^{*2}\tau^2+2\alpha^*\tau-1)}\neq 0,
\label{condeq20}
\end{equation}
(where ${\cal G}$ is as in (\ref{cond3})), and that the first Lyapunov coefficient $K_1^*\equiv K_1(\alpha^*,\beta^*)$ computed in Appendix A is non-zero.  Then (\ref{SIS}) has a degenerate Hopf bifurcation point at the endemic equilibrium $\bar{y}$ when $(R_0,p)=(R_0^*,p^*)$, and the normal form for the Hopf bifurcation diagrams near this point is
\[
r(\varepsilon (\lambda^2+\eta)+r^2)=0
\]
where $\lambda=R_0-R_0^*$, $\varepsilon=\mbox{\rm sgn}(\sigma_4^*/K_1^*)$ and 
${\displaystyle\eta=\frac{\sigma_1^* (p-p^*)}{|K_1^*|\mbox{\rm sgn}(\sigma_4^*)}}$.  Furthermore, for $(R_0,p)$ close enough to $(R_0^*,p^*)$, this Hopf bifurcation diagram will exhibit an endemic bubble if $\varepsilon=+1$ and $\eta<0$.  In the case $\varepsilon=+1$ and $\eta<0$, the width of the endemic bubble (in $R_0$ space) is approximately equal to the width of the region between the two zeros of
\[
\sigma_4(R_0-R_0^*)^2+\sigma_1 (p-p^*),
\]
i.e.
\[
\tilde{R}_0-\bar{R}_0\approx 2\sqrt{\left|\frac{\sigma_1(p-p^*)}{\sigma_4}\right|}
\]
\label{thm42}
\end{thm}

\subsection{Response function $h(y,p)=\frac{1}{1+py}$}

This is one of the special cases which was studied in \cite{bubbles} and for which endemic bubbles were observed.  For purposes of comparing our results to those of \cite{bubbles}, we will assume as they do that $\tau=10$.

\begin{figure}[htpb]
\begin{center}
\includegraphics[width=3.0in]{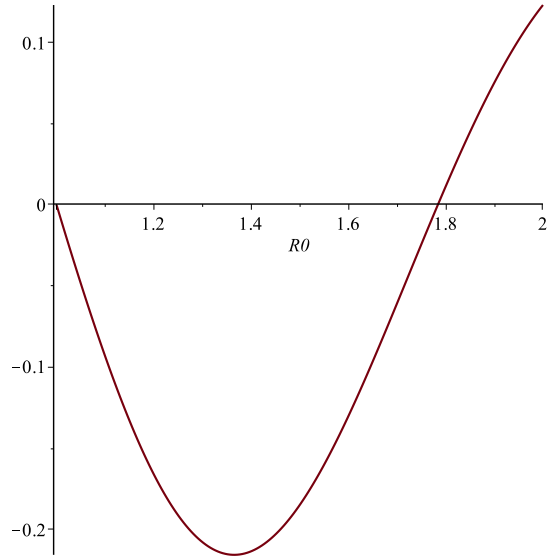}
\caption{Graph of the left-hand side of (\ref{condeq5}).}
\label{fig7}
\end{center}
\end{figure}

The endemic equilibrium is
\[
\bar{y}(R_0,p)=\frac{R_0-1}{p+R_0},\,\,\,\,R_0>1
\]
and introducing $x=y-y^*$ transforms (\ref{SIS}) into
\begin{equation}
\dot{x}(t)=\alpha(R_0,p)x(t)+\beta(R_0,p)x(t-\tau)+O(|x|^2),
\label{translate1}
\end{equation}
where 
\[
\alpha(R_0,p)=\frac{1-R_0}{1+p},\,\,\,\,\,\beta(R_0,p)=\frac{(1-R_0)p}{R_0(1+p)}.
\]
Computing the quantity $\beta\alpha_{R_0}(1-\alpha\tau)+\beta_{R_0}(\tau\beta^2-\alpha)$ for $\tau=10$ gives
\[
{\frac { \left( 1-{R_0} \right) ^{2}p \left( 10\,{{R_0}}^{3}+{{
R_0}}^{2}p+{{R_0}}^{2}-10\,{p}^{2} \right) }{{{R_0}}^{4}
 \left( 1+p \right) ^{3}}}.
\]
So using (\ref{condeq2}) and solving for $p$, we have
\[
p=p^*=\frac{1}{20}\, \left( {R_0^*}+\sqrt {{{R_0^*}}^{2}+400\,{R_0^*}+40}
 \right) {R_0^*},
 \]
 where $R_0^*$ is a root of (\ref{condeq1}), which is equivalent to
 \begin{equation}
 \cos \left( 10\,\omega^* \right) \left( R_0^*+\sqrt {{{R_0^*}}^{2}+400\,{R_0^*}+40}\right)
+20=0
\label{condeq5}
\end{equation}
where
\[
\omega^{*2}=2\,{\frac { \left( 1-{R_0^*} \right) ^{2} \left( {R_0^*}\,\sqrt {{{
R_0^*}}^{2}+400\,{R_0^*}+40}+{{R_0^*}}^{2}+200\,{R_0^*}-180
 \right) }{ \left( {R_0^*}\,\sqrt {{{R_0^*}}^{2}+400\,{R_0^*}+40}+{
{R_0^*}}^{2}+20 \right) ^{2}}}.
\]
A plot of the left-hand side of (\ref{condeq5}) as a function of $R_0$ is given in Figure \ref{fig7}.  We numerically compute the value for the root
\[
R_0^*\approx 1.784
\]
which yields
\[
p^*\approx 2.613
\]
and
\[
\alpha^*\approx -0.217,\,\,\,\beta^*\approx -0.318,\,\,\,\,\omega^*\approx 0.232.
\]
We compute $\sigma_1^*$ and $\sigma_4^*$ as in (\ref{condeq10}) and (\ref{condeq20}) and get
\[
\sigma_1^*\approx 0.021,\,\,\,\,\,\,\,\,\,\,\,\,\,\,\sigma_4^*\approx -0.037.
\]

To compute the Lyapunov coefficient $K_1^*$, we will need to compute explicitly the quadratic and cubic terms in (\ref{translate1}).  
A straightforward computation gives the following expression for these quadratic and cubic terms:
\[
\begin{array}{l}
-\left(\frac{p^*+R_0^*}{1+p^*}\right)x(t)^2+
\left(\frac {p^*(R_0^*-2-p^*)(p^*+R_0^*)}{R_0^*(1+p^*)^2}\right)
x(t)x(t-\tau)
+\left(\frac{(p^*+R_0^*)(R_0^*-1)p^{*2}}{R_0^{*2}(1+p^*)^2}\right)
x(t-\tau)^2\\[0.2in]
+\left(\frac{(p^*+R_0^*)^2p^*}{R_0^*(1+p^*)^2}\right)
x(t)^2x(t-\tau)-
\left(\frac{p^{*2}(p^*+R_0^*)^2(R_0^*-2-p^*)}{R_0^{*2}(1+p^*)^3}\right)
x(t)x(t-\tau)^2\\[0.2in]
+\left(\frac{(1-R_0^*)(p^*+R_0^*)^2p^{*3}}{R_0^{*3}(1+p^*)^3}\right)x(t-\tau)^3.
 \end{array}
\]
The formula for $K_1$ in the Appendix now gives
\[
K_1^*\approx -1.006
\]
so that 
\[
\varepsilon=\mbox{\rm sgn}(\sigma_4^*/K_1^*) =+1,\,\,\,\,\,\,\,\,\,\eta\approx -0.021\, (p-p^*).
\]
Based on Theorem \ref{thm42}, we therefore predict the following for $p$ near $p^*$ when $h(y,p)=1/(1+yp)$:
\begin{enumerate}
\item[(i)] Suppose $p<p^*$. Then varying $R_0$ near $R_0^*$ in (\ref{SIS}) will not generate Hopf bifurcation from the endemic equilibrium, hence there will be no endemic bubble.
\item[(ii)] Suppose $p>p^*$.  Then varying $R_0$ near $R_0^*$ in (\ref{SIS}) will generate an endemic bubble via Hopf bifurcation at $R_0=\overline{R}_0$ followed by a reverse Hopf bifurcation at $R_0=\tilde{R}_0$, with $\overline{R}_0$ and $\tilde{R}_0$ near $R_0^*$, and
\[
\tilde{R}_0-\overline{R}_0\approx 2\sqrt{-\frac{\sigma_1^*}{\sigma_4^*}\,(p-p^*)}\approx 1.614\sqrt{p-p^*}.
\]
\end{enumerate}
\begin{figure}[htpb]
\begin{center}
\includegraphics[width=2.5in]{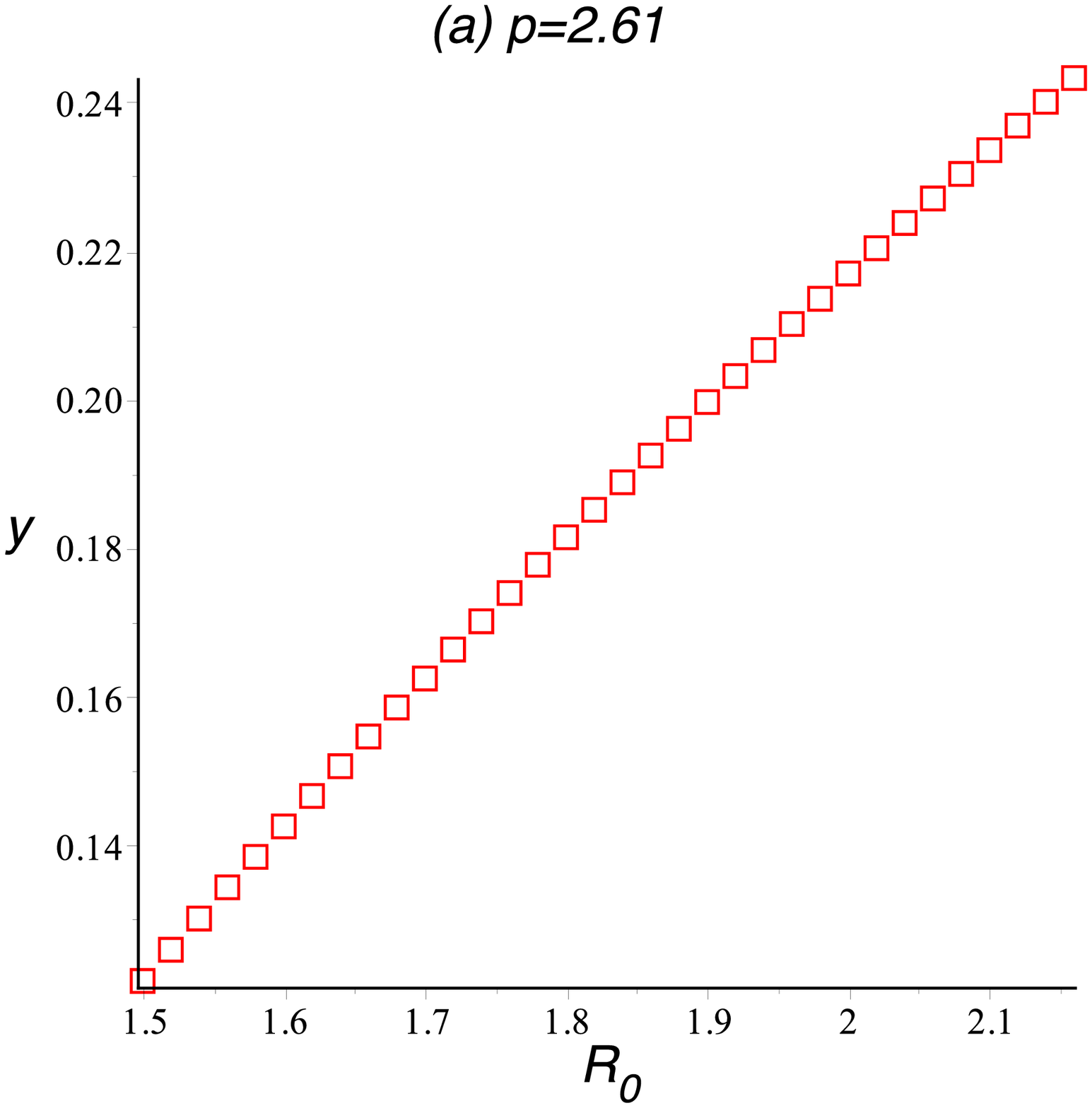}
\includegraphics[width=2.5in]{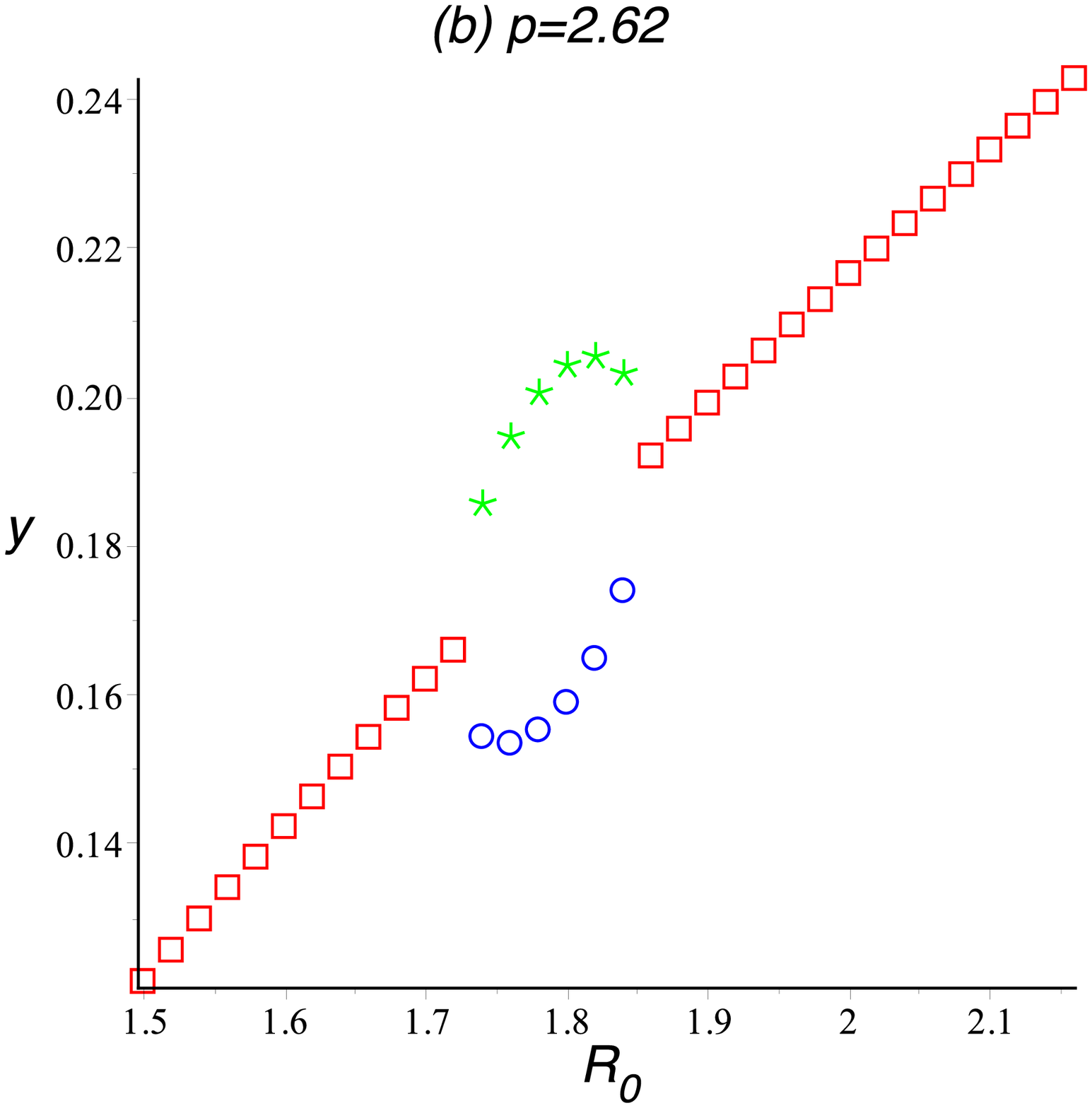}

\includegraphics[width=2.5in]{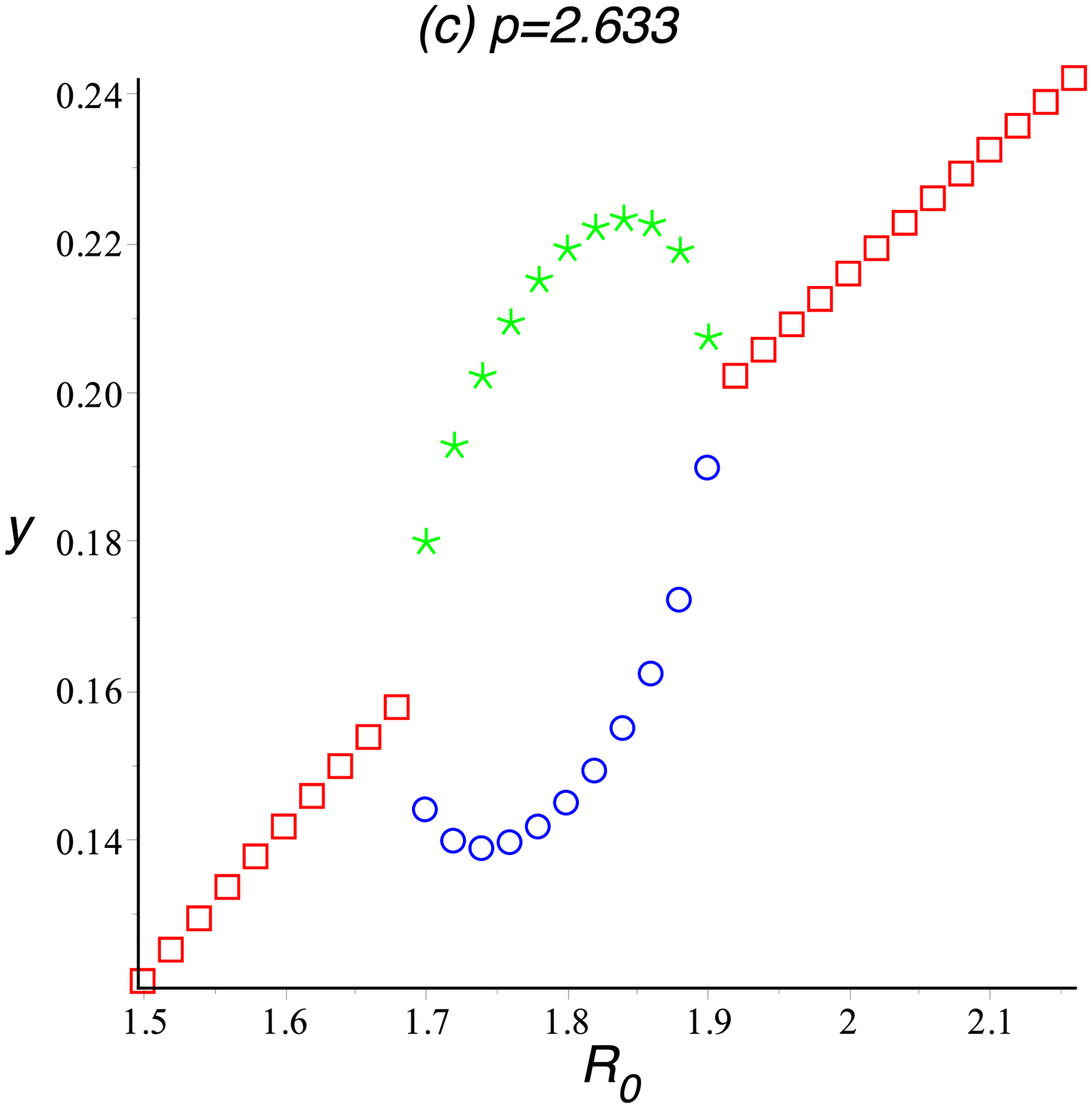}
\includegraphics[width=2.5in]{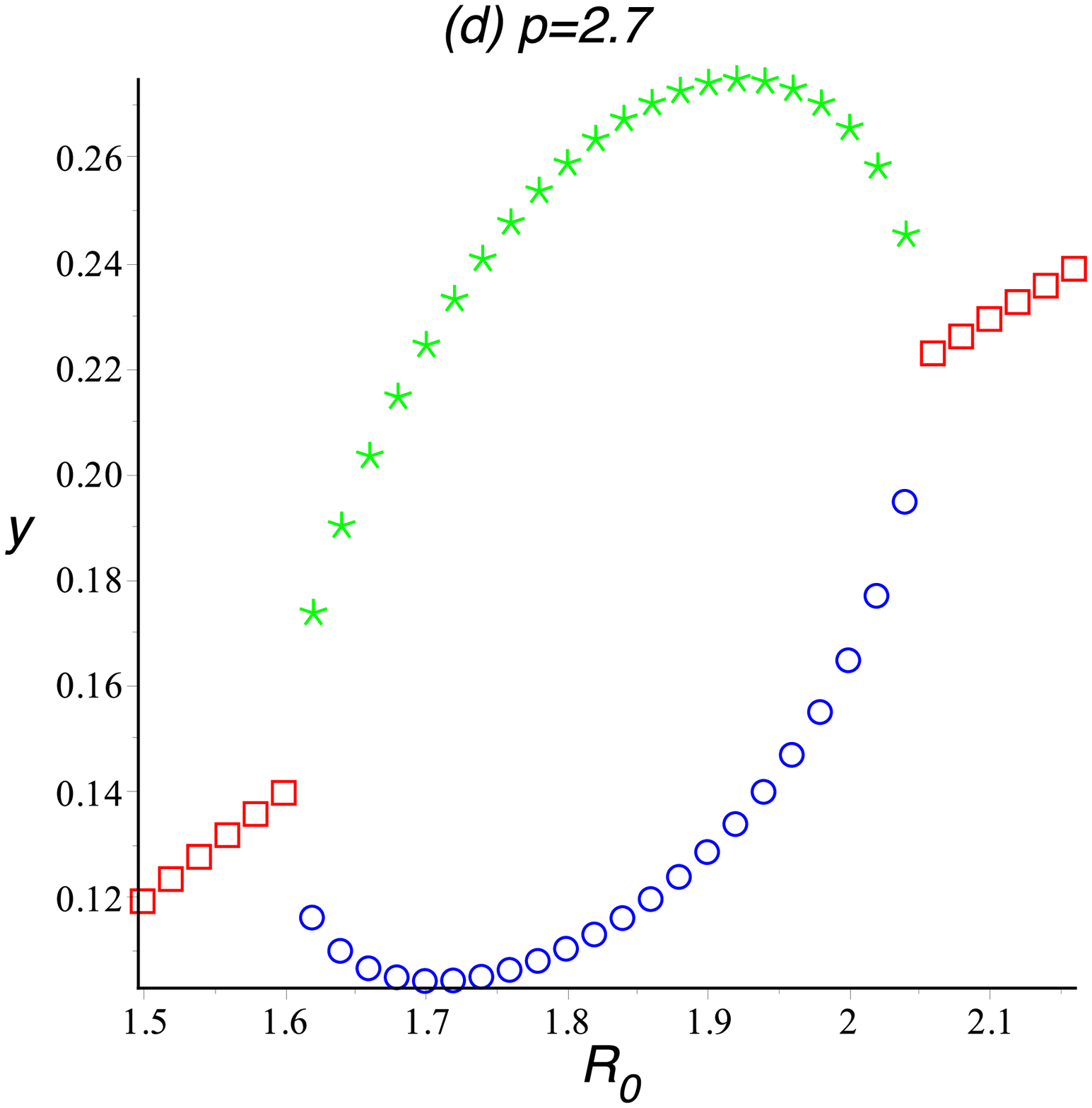}
\caption{Numerically obtained bifurcation diagrams for (\ref{SIS}) with $h(y,p)=1/(1+py)$ for the following values of $p$: (a) $p=2.61<p^*$, (b) $p=2.62>p^*$, (c) $p=2.633>p^*$, (d) $p=2.7>p^*$.  The red squares indicate that the simulation has settled to the endemic equilibrium.  The green asterisks and blue circles designate respectively the maxima and minima of the steady oscillations which occur for values of $R_0$ inside the endemic bubble. Recall that the theoretically computed value of $(R_0^*,p^*)$ for the degenerate Hopf bifurcation is $(R_0^*,p^*)\approx (1.784,2.613)$.}
\label{fig8}
\end{center}
\end{figure}
\begin{figure}[htpb]
\begin{center}
\includegraphics[width=2.5in]{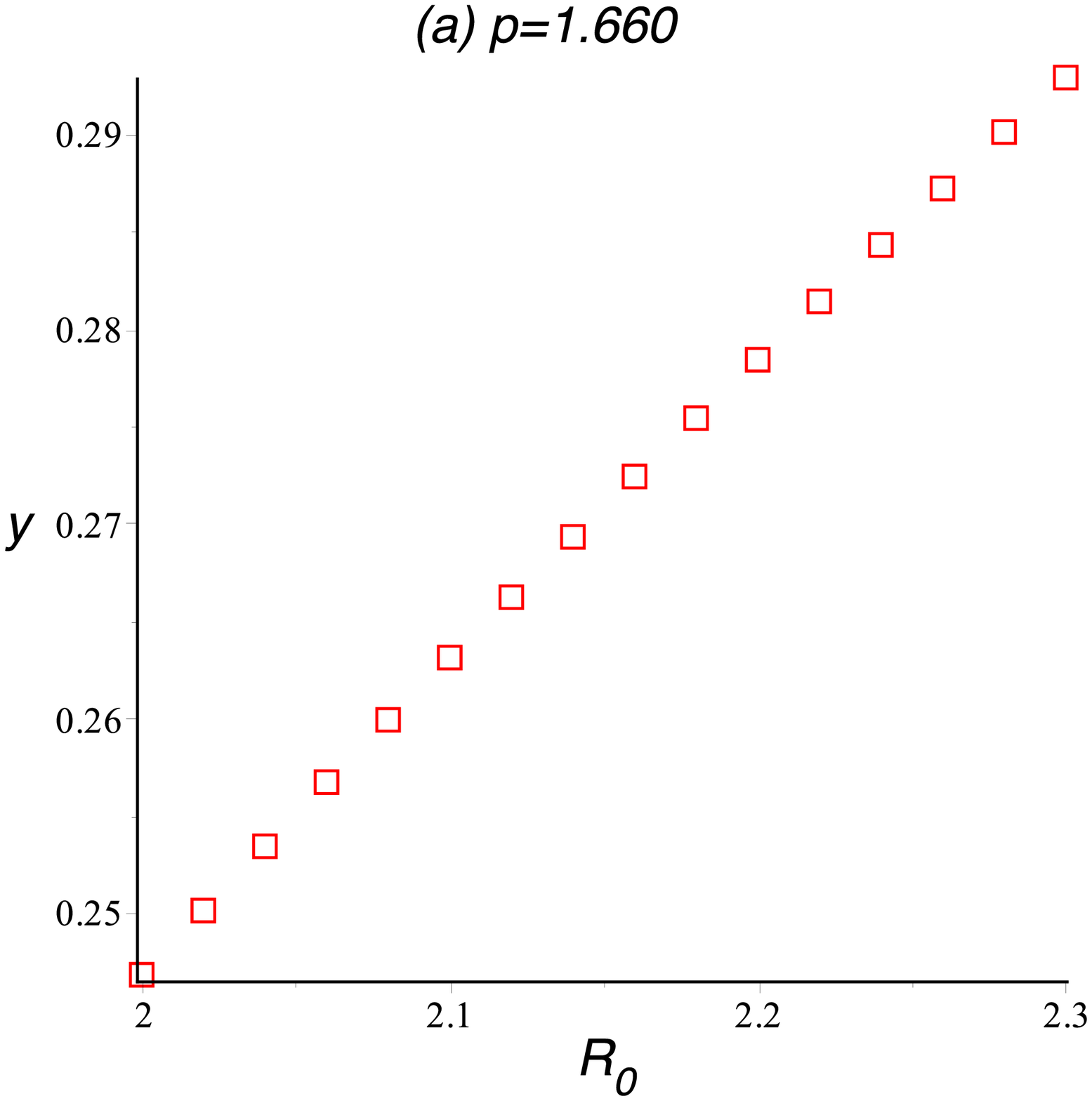}
\includegraphics[width=2.5in]{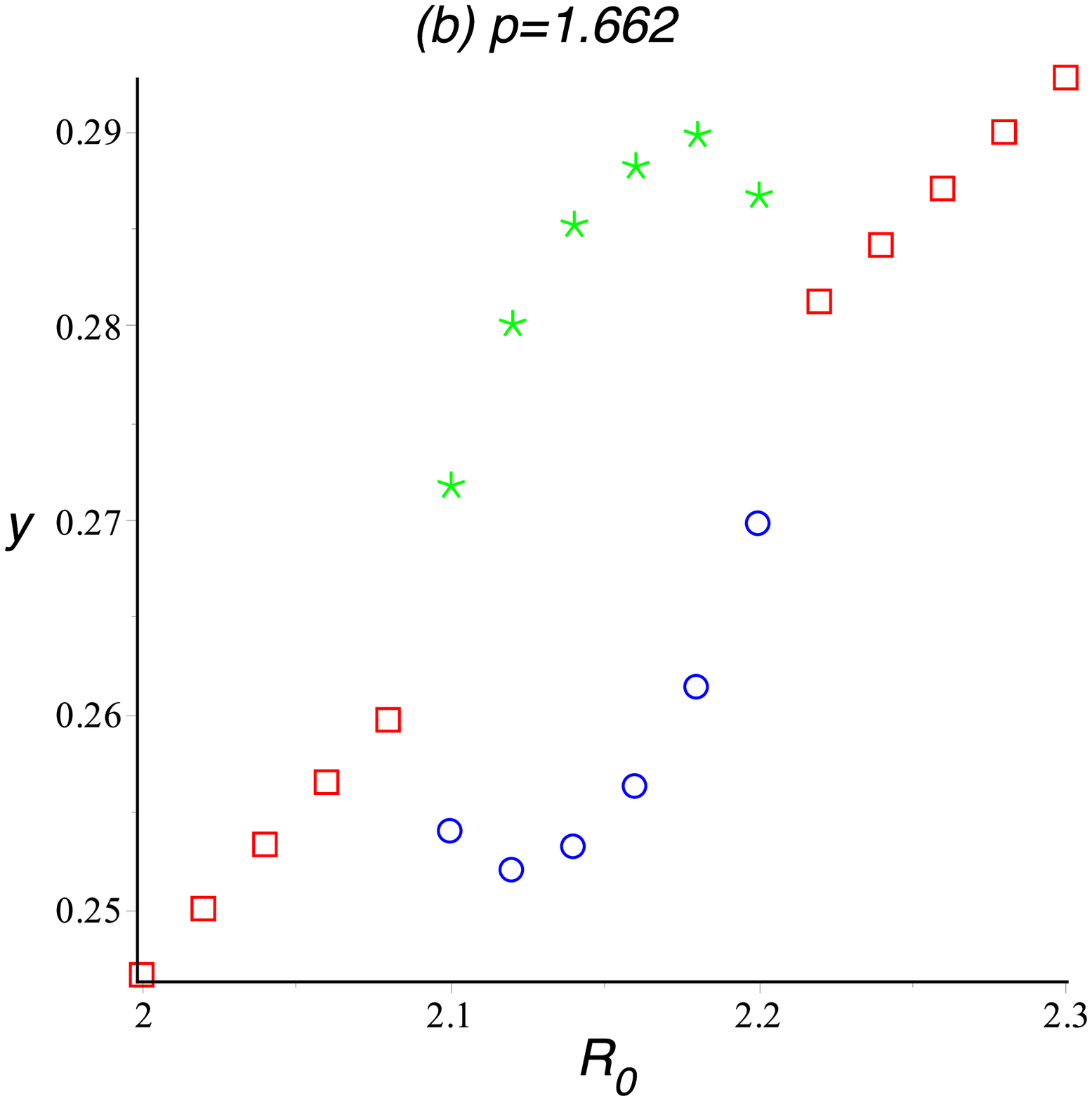}
\caption{Numerically obtained bifurcation diagrams for (\ref{SIS}) with $h(y,p)=e^{-py}$ for the following values of $p$: (a) $p=1.660<p^*$, (b) $p=1.662>p^*$.  The red squares indicate that the simulation has settled to the endemic equilibrium.  The green asterisks and blue circles designate respectively the maxima and minima of the steady oscillations which occur for values of $R_0$ inside the endemic bubble. Recall that the theoretically computed value of $(R_0^*,p^*)$ for the degenerate Hopf bifurcation is $(R_0^*,p^*)\approx (2.1474,1.6617)$.}
\label{fig9}
\end{center}
\end{figure}
\begin{figure}[htpb]
\begin{center}
\includegraphics[width=5in]{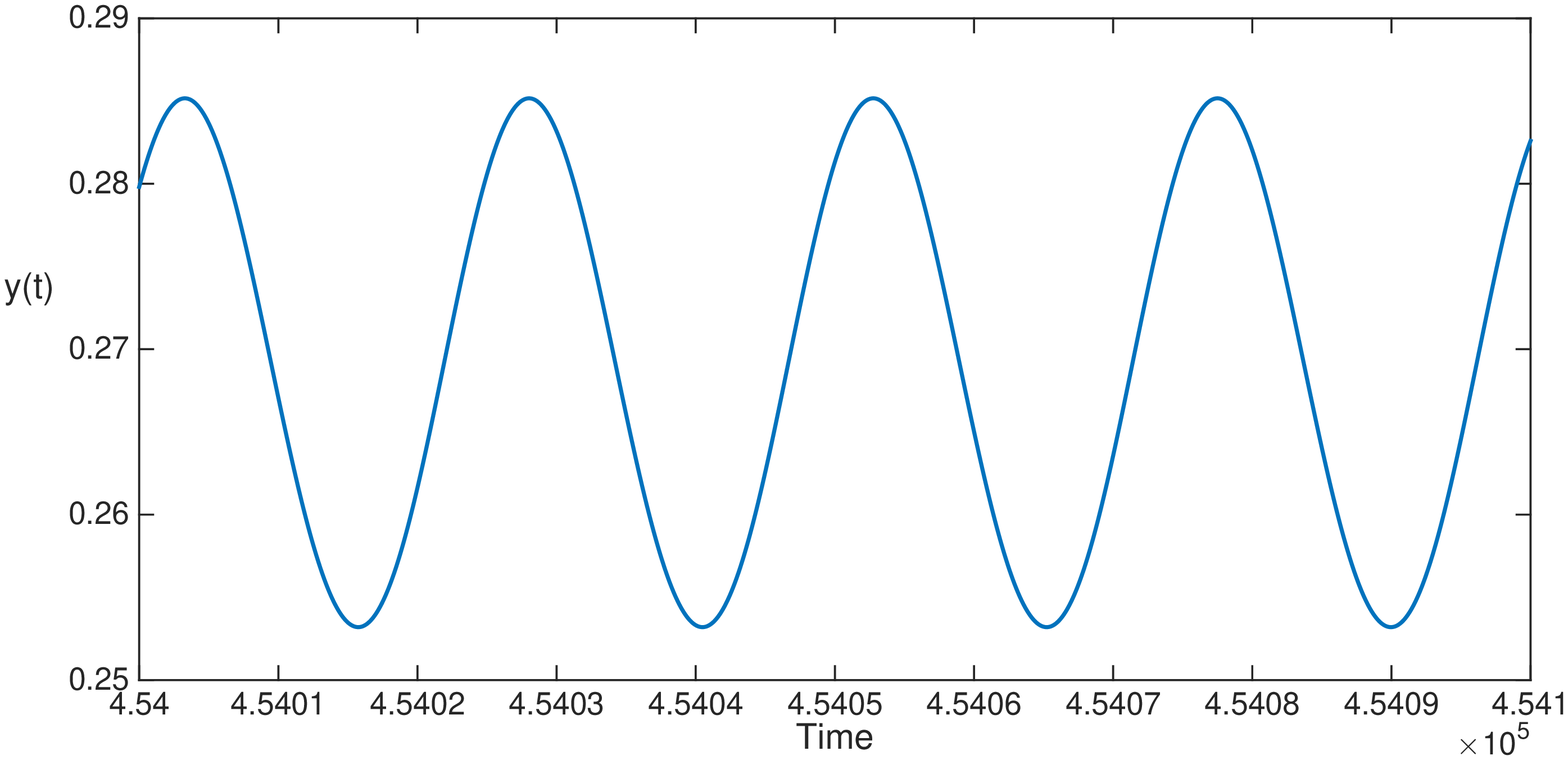}
\caption{Numerically obtained periodic solution $y(t)$ to (\ref{SIS}) with $h(y,p)=e^{-py}$ after transients have died out. The parameter values are $p=1.662$ and $R_0=2.14$, which is inside the endemic bubble.  The period is approximately 25, which is close to the theoretically predicted value $2\pi/\omega^*\approx 24.7$.}
\label{fig10}
\end{center}
\end{figure}

In Figure \ref{fig8}, we illustrate the results of several integrations of (\ref{SIS}) with $h(y,p)=1/(1+py)$ for values of $(R_0,p)$ near $(R_0^*,p^*)$.  The resulting bifurcation diagrams show the emergence of the endemic bubble.

\subsection{Response function $h(y,p)=e^{-py}$}

This case was also studied in \cite{bubbles}, although the bifurcation diagrams were not given in that paper.  In this case, we can not solve in closed form the equation for the endemic equilibrium (\ref{eqeq3})
\begin{equation}
e^{p\bar{y}(R_0,p)}=R_0(1-\bar{y}(R_0,p)).
\label{neweqeq3}
\end{equation}
Using (\ref{neweqeq3}), equations (\ref{linalpha}) and (\ref{linbeta}) become
\[
\alpha(R_0,p)=\frac{-\bar{y}(R_0,p)}{1-\bar{y}(R_0,p)}
\]
and
\[
\beta(R_0,p)=-p\bar{y}(R_0,p).
\]
Implicit differentiation of (\ref{neweqeq3}) gives
\[
\bar{y}_{R_0}(R_0,p)=\frac{1-\bar{y}(R_0,p)}{R_0(1+p(1-\bar{y}(R_0,p)))},
\]
from which we may now compute
\[
\alpha_{R_0}(R_0,p)=-{\frac {1}{ \left( 1-\bar{y}(R_0,p) \right) {R_0}\, \left( 1+p(1-\bar{y}(R_0,p)) \right) }}
\]
and
\[
\beta_{R_0}(R_0,p)=-{\frac {p \left( 1-\bar{y}(R_0,p) \right) }{{R_0}\, \left( 1+p(1-\bar{y}(R_0,p)) \right) }}.
\]
For $\tau=10$, solving (\ref{condeq2}) becomes equivalent to solving
\[
-10(1-\bar{y}(R_0,p))^3\,p^2+11-\bar{y}(R_0,p)=0,
\]
or
\begin{equation}
p=\frac{1}{10}\,\frac{\sqrt{10 (1-\bar{y}(R_0,p))(11-\bar{y}(R_0,p))}}{(1-\bar{y}(R_0,p))^2}.
\label{pcomp}
\end{equation}
Using this expression for $p$, we may solve (numerically) the equation
\[
\alpha+\beta\,\cos(10\sqrt{\beta^2-\alpha^2})=0
\]
for $\bar{y}$, and the result is
\[
\bar{y}\approx 0.2703
\]
which, when substituted into (\ref{pcomp}) gives
\[
p^*\approx 1.6617
\]
and then (\ref{neweqeq3}) gives
\[
R_0^*\approx 2.1474.
\]
We then compute
\[
\alpha^*\approx -0.3704,\,\,\,\,\,\,\,\,\beta^*\approx -0.4491,\,\,\,\,\,\,\,\,\omega^*\approx 0.2540.
\]
Using these values and further implicit differentiations of (\ref{neweqeq3}) (we omit the details) yields
\[
\sigma_1^*\approx 0.0503,\,\,\,\,\,\,\,\,\,\,\,\,\sigma_4^*\approx -0.0190.
\]
Finally, the first Lyapunov coefficient has the value
\[
K_1^*\approx -0.4906
\]
so that
\[
\varepsilon=\mbox{\rm sgn}(\sigma_4^*/K_1^*)=+1,\,\,\,\,\,\,\,\eta\approx -0.1025\,(p-p^*).
\]
Based on Theorem \ref{thm42}, we therefore predict the following for $p$ near $p^*$ when $h(y,p)=e^{-py}$:
\begin{enumerate}
\item[(i)] Suppose $p<p^*$. Then varying $R_0$ near $R_0^*$ in (\ref{SIS}) will not generate Hopf bifurcation from the endemic equilibrium, hence there will be no endemic bubble.
\item[(ii)] Suppose $p>p^*$.  Then varying $R_0$ near $R_0^*$ in (\ref{SIS}) will generate an endemic bubble via Hopf bifurcation at $R_0=\overline{R}_0$ followed by a reverse Hopf bifurcation at $R_0=\tilde{R}_0$, with $\overline{R}_0$ and $\tilde{R}_0$ near $R_0^*$, and
\[
\tilde{R}_0-\overline{R}_0\approx 2\sqrt{-\frac{\sigma_1^*}{\sigma_4^*}\,(p-p^*)}\approx 4.486\sqrt{p-p^*}.
\]
\end{enumerate}
Resuts of numerical simulations for this case are illustrated in Figures \ref{fig9} and \ref{fig10}.

\Section{Conclusions}

In this paper, we have performed a detailed theoretical analysis of a degenerate Hopf bifurcation in parametrized families of RFDEs where the degeneracy arises from a violation of the eigenvalue crossing condition.  Although our detailed computations have been performed for scalar RFDEs, there are no complications other than those involved with cumbersome notation in extending these results to systems of RFDEs.   We have paid particular attention to the cases where the RFDE is a delay differential equation of the form (\ref{prot}), since these cases are quite ubiquitous in the literature, and many important phenomena in nature are modelled using such equations.  Since we give precise conditions on the parameters of (\ref{prot}) to characterize the degenerate Hopf bifurcation, we expect that our paper will be a valuable contribution to many researchers who model phenomena using delay differential equations.

As an application of our results, we have considered the SIS model (\ref{SIS}) which was studied in \cite{bubbles}.  We have shown that the phenomenon of endemic bubbles which had been reported in \cite{bubbles} originates in a degenerate Hopf bifurcation such as the one studied herein.  

It is worth noting that the formula for the Lyapunov coefficient in the appendix includes the parameters $\alpha^*$, $\beta^*$ and $\tau$.  It is conceivable that by varying $\tau$, we could achieve an even higher-order degeneracy where both the crossing condition and the first Lyapunov condition are violated.  This is work in progress.

\appendix

\Section{First Lyapunov coefficient}

For the delay-differential equation (\ref{prot}), formulae (\ref{Kdef}), (\ref{Kdef2}) and (\ref{Fdef}) yield the Lyapunov coefficient $K_1$ in terms of the coefficients $f_{(j,k)}$ in (\ref{Fdef}) and $\alpha^*$, $\beta^*$, $\omega^*$ and $\tau$.  The following result was computed using the symbolic computation software package Maple:

\[
\begin{array}{l}
K_1=\frac{1}{(1-\alpha^*\tau)^2+\omega^{*2}}\left[3(1-\alpha^*\tau)f_{(3,0)}
+\left({\frac {3\,{\alpha^*}^{2}\tau-{\alpha^*}^{2}+{\beta^*}^{2}-3\,\alpha^*}{\beta^*}}\right)f_{(2,1)}\right.\\[0.2in]
-\left({\frac {2\,{\alpha^*}^{3}\tau+\alpha^*\,{\beta^*}^{2}\tau-2\,{\alpha^*}^{3}+2
\,\alpha^*\,{\beta^*}^{2}-2\,{\alpha^*}^{2}-{\beta^*}^{2}}{{\beta^*}^{2}}}\right)f_{(1,2)}
+3\left(\,{\frac {{\alpha^*}^{2}\tau-{\alpha^*}^{2}+{\beta^*}^{2}-\alpha^*}{\beta^*}}\right)f_{(0,3)}\\[0.2in]
+2\left(\,{\frac {6\,{\alpha^*}^{2}\tau-9\,\alpha^*\,\beta^*\,\tau-2\,{\alpha^*}^{2}+
2\,{\beta^*}^{2}-6\,\alpha^*+9\,\beta^*}{ \left( \alpha^*+\beta^* \right) 
 \left( 4\,\alpha^*-5\,\beta^* \right) }}\right)f_{(2,0)}^2\\[0.2in]
-\left({\frac {18\,{\alpha^*}^{3}\tau-33\,{\alpha^*}^{2}\beta^*\,\tau+9\,\alpha^*\,{
\beta^*}^{2}\tau-10\,{\alpha^*}^{3}+7\,{\alpha^*}^{2}\beta^*+10\,\alpha^*\,{
\beta^*}^{2}-7\,{\beta^*}^{3}-18\,{\alpha^*}^{2}+33\,\alpha^*\,\beta^*-9\,{\beta^*
}^{2}}{ \left( \alpha^*+\beta^* \right)  \left( 4\,\alpha^*-5\,\beta^*
 \right) \beta^*}}\right)f_{(2,0)}f_{(1,1)}\\[0.2in]
-2\left(\,{\frac { \left( \alpha^*-\beta^* \right)  \left( 6\,{\alpha^*}^{2}\tau-9
\,\alpha^*\,\beta^*\,\tau-6\,{\alpha^*}^{2}+\alpha^*\,\beta^*+7\,{\beta^*}^{2}-6\,
\alpha^*+9\,\beta^* \right) }{ \left( \alpha^*+\beta^* \right)  \left( 4\,
\alpha^*-5\,\beta^* \right) \beta^*}}\right)f_{(2,0)}f_{(0,2)}\\[0.2in]
+
\left({\frac { \left( \alpha^*-\beta^* \right)  \left( 4\,{\alpha^*}^{3}\tau-10\,{
\alpha^*}^{2}\beta^*\,\tau+\alpha^*\,{\beta^*}^{2}\tau-4\,{\alpha^*}^{3}+2\,{
\alpha^*}^{2}\beta^*+3\,\alpha^*\,{\beta^*}^{2}-3\,{\beta^*}^{3}-4\,{\alpha^*}^{2}
+10\,\alpha^*\,\beta^*-{\beta^*}^{2} \right) }{{\beta^*}^{2} \left( \alpha^*+
\beta^* \right)  \left( 4\,\alpha^*-5\,\beta^* \right) }}\right)f_{(1,1)}^2\\[0.2in]

+\left({\frac {8\,\tau\,{\alpha^*}^{5}+8\,{\alpha^*}^{4}\beta^*\,\tau-32\,{\alpha^*}^
{3}{\beta^*}^{2}\tau+19\,{\alpha^*}^{2}{\beta^*}^{3}\tau-9\,\alpha^*\,{\beta^*}^
{4}\tau-8\,{\alpha^*}^{5}-8\,{\alpha^*}^{4}\beta^*+36\,{\alpha^*}^{3}{\beta^*}^{
2}}{{\beta^*}^{3} \left( \alpha^*+\beta^* \right)  \left( 4\,\alpha^*-5\,\beta^*
 \right) }}\right.\\[0.2in]
 
\,\,\,\,\,\,\,\,\,\,\,\,\,\,\,\,\,\,\,\,\,\,\,\,\,\,\,\,\,\,\,\,\,\,\,\,\,\,\,\,\,\,\,\,\,\,\,\,\,\,\,\,\,\,\,\,\,\,\,\,\,\, +\left.{\frac {{\alpha^*}^{2}{\beta^*}^{3}-28\,\alpha^*\,{\beta^*}^{4}+7\,{\beta^*}^{5}
-8\,{\alpha^*}^{4}-8\,{\alpha^*}^{3}\beta^*+32\,{\alpha^*}^{2}{\beta^*}^{2}-19\,
\alpha^*\,{\beta^*}^{3}+9\,{\beta^*}^{4}}{{\beta^*}^{3} \left( \alpha^*+\beta^*
 \right)  \left( 4\,\alpha^*-5\,\beta^* \right) }}\right)f_{(1,1)}f_{(0,2)}\\[0.2in]

-2\left(\frac{4\alpha^{*4}\tau+4\alpha^{*3}\beta^*\tau-13\alpha^{*2}\beta^{*2}\tau+2\alpha^*\beta^{*3}\tau-4\alpha^{*4}-4\alpha^{*3}\beta^*+15\alpha^{*2}\beta^{*2}+4\alpha^*\beta^{*3}}{\beta^{*2}(\alpha^*+\beta^*)(4\alpha^*-5\beta^*)}\right.\\[0.2in]

\,\,\,\,\,\,\,\,\,\,\,\,\,\,\,\,\,\,\,\,\,\,\,\,\,\,\,\,\,\,\,\,\,\,\,\,\,\,\,\,\,\,\,\,\,\,\,\,\,\,\,\,\,\,\,\,\,\,\,\,\,\,\,\,\,\,\,\,\,\,\,\,\,\,\,\,\,\,\,+
\left.\left.\frac{-11\beta^{*4}-4\alpha^{*3}-4\alpha^{*2}\beta^*+13\alpha^*\beta^{*2}-2\beta^{*3}}{\beta^{*2}(\alpha^*+\beta^*)(4\alpha^*-5\beta^*)}\right)\,f_{(0,2)}^2\right].
%
%
\end{array}
\]

It follows that in order for $K_1$ to be well-defined, we need the non-degeneracy conditions $\alpha^*+\beta^*\neq 0$ and $4\alpha^*-5\beta^*\neq 0$.  For a generic delay-differential equation of the form (\ref{prot}), this coefficient $K_1$ will be non-zero.

\vspace*{0.25in}
\noindent
{\Large\bf Acknowledgments}

\vspace*{0.2in}
This research is partly supported by the
Natural Sciences and Engineering Research Council of Canada in the
form of a Discovery Grant.  The author is grateful to the reviewers for carefully reading the manuscript and
for providing valuable suggestions that have considerably improved the paper.

\end{document}